\documentclass[pdflatex,11pt]{article}

\usepackage{a4wide}
\usepackage{amssymb}
\usepackage{latexsym,graphics,tikz}
\usepackage{amsfonts,amssymb,amsmath,mathtools}
\usepackage[space]{grffile}
\usepackage{epstopdf}
\usepackage{mathrsfs}
\usepackage{graphicx}
\usepackage{subfig}
\usepackage[figuresright]{rotating}

\graphicspath{{../MATLAB_code/}}

\newcommand{\ZZ}{\mathbb{Z}}

\newcommand{\IH}{\mathbb{H}}

\newcommand{\dx}{\Delta x}
\newcommand{\dt}{\Delta t}
\newcommand{\zni}{z_{n}^{i}}

\newcommand{\half}{\mbox{$\frac{1}{2}$}}
\newcommand{\ti}{\mathtt{i}}

\DeclareMathOperator{\sech}{sech}

\begin{document}

%\begin{frontmatter}


%% \author{A.~Bhatt\corref{cor1}\fnref{label2}}

\begin{center}
{\Large\bf Exponential Integrators Preserving Local Conservation Laws of PDEs with 
Time-Dependent Damping/Driving Forces}
%with \\Application to Damped-Driven NLS Equations}

\renewcommand{\thefootnote}{\fnsymbol{footnote}}

\vspace{4ex}
{\sc Ashish Bhatt\footnotemark[1] 
and Brian E.~Moore\footnotemark[2]}

\footnotetext[1]{Department of Mathematics, University of Stuttgart,
($\mathtt{ashish.bhatt@mathematik.uni-stuttgart.de}$).}

% \footnotetext[2]{School of Science, Chang'an University,
% ($\mathtt{yuyue.qin@chd.edu.cn}$).}

\footnotetext[2]{Department of Mathematics, University of Central Florida,
($\mathtt{brian.moore@ucf.edu}$).}

\renewcommand{\thefootnote}{\arabic{footnote}}
\end{center}


%\author[Bhatt]{Ashish Bhatt} 
%\author[Moore]{Brian E.~Moore}
%\maketitle
%\address[Bhatt]{Department of Mathematics, University of Stuttgart, ($\mathtt{ashish.bhatt@mathematik.uni-stuttgart.de}$)}
%\address[Moore]{Department of Mathematics, University of Central Florida, ($\mathtt{brian.moore@ucf.edu}$)}

\begin{abstract}
Structure-preserving algorithms for solving conservative PDEs with added linear dissipation are 
generalized to systems with time dependent damping/driving terms. 
This study is motivated by several PDE models of physical phenomena,  
such as Korteweg-de Vries, Klein-Gordon, Schr\"{o}dinger, and Camassa-Holm equations, 
all with damping/driving terms and time-dependent coefficients.
Since key features of the PDEs under consideration are described by local conservation laws,
which are independent of the boundary conditions,
the proposed (second-order in time) discretizations are developed with the intent of preserving those 
local conservation laws. 
The methods are respectively applied to a damped-driven nonlinear Schr\"{o}dinger equation 
and a damped Camassa-Holm equation.  
Numerical experiments illustrate the structure-preserving properties 
of the methods, as well as favorable results over other competitive schemes.
\end{abstract}

\vspace{4ex}

\vspace{2ex}

%\end{frontmatter}

% \linenumbers


%%%%%%%%%%%%%%%%%%%%%%%%%%%%%%%%%%%%%%%%%%%%%%%%%%%%%%%%%%%%%%%%%%%%%%%%
\section{Introduction} \label{sec:intro}
%%%%%%%%%%%%%%%%%%%%%%%%%%%%%%%%%%%%%%%%%%%%%%%%%%%%%%%%%%%%%%%%%
When modeling physical phenomena, it is important to maintain as many properties of the physical system as possible.
Numerical models that fail to maintain aspects of the physical behavior are actually modeling something other than the 
intended physics.  This idea has been instrumental in the development of structure-preserving algorithms for solving 
differential equations that model conservative dynamics.
Many successful time-stepping schemes for solving conservative differential equations are constructed to 
ensure that aspects of the qualitative solution behavior are preserved. 
However, many physical systems are commonly effected by external forces or by the dissipative effects of friction, and 
methods that preserved the structure in standard conservative systems fail to preserve the structure when damping/driving forces are included in the model.  

The need for numerical methods that preserve dissipative properties of physical systems has resulted in the development of 
some useful discretization techniques.
In some cases, guaranteeing monotonic dissipation of a desired quantity (numerically) may be sufficient (cf. \cite{Cetal12}).
In this work we focus on a special form of linear damping that can be preserved more precisely. 
Numerical methods that preserve the conformal symplectic structure of conformal Hamiltonian systems \cite{MP01}
are naturally known as conformal symplectic methods.  Such methods were first constructed for ODEs
using splitting techniques \cite{McLQ01,McLQ02}, by solving the dissipative part exactly and 
the conservative part with a symplectic method, and then composing the flow maps.  
In subsequent work, various second order conformal symplectic methods were proposed \cite{BFM15,MS11}, 
which included linear stability analysis, backward error analysis, additional results on 
structure-preservation, and favorable numerical results for both ODE and PDE examples. 

In fact, there are several works that use these ideas 
to construct structure preserving methods for PDEs.  
Conformal symplectic Euler methods were used to solve damped nonlinear wave equations \cite{Moore09}. 
Then, the ideas were extended to general linearly damped multi-symplectic PDEs, 
with application to damped semi-linear waves equations and damped NLS equations, 
where midpoint-type discretizations 
(similar to the Preissmann box scheme and midpoint discrete gradient methods) 
were shown to preserve various dissipative properties of the governing equations \cite{MNS13}.
Other similar methods have been developed based on 
an even-odd Strang splitting for a modfied (non-diffusive) Burger's equation \cite{BFM15}, 
the average vector field method  for a damped NLS equation \cite{Fu16}, and an 
explicit fourth-order Nystr\"{o}m method with composition techniques for damped acoustic wave equations \cite{CZW17}.

More recently, conformal symplectic integration has been generalized in a few different ways.
All of the aforementioned works on methods for PDEs have hinged on conformal symplectic time discretizations,
but the ideas have been fully generalized to PDEs, including conformal symplectic spatial discretizations, 
with several example model problems \cite{M16}.
In the context of conformal Hamiltonian ODEs, general classes of exponential Runge-Kutta methods
and partitioned Runge-Kutta methods provide higher order methods, where
sufficient conditions on the coefficient functions guarantee 
preservation of conformal symplecticity or other quadratic conformal invariants \cite{BM17}. 
The framework proposed there also provides straightforward approaches to constructing 
conformal symplectic methods, using a Lawson-type transformation.  
As shown in \cite{BM17}, many of the conformal symplectic exponential Runge-Kutta methods are 
applicable to differential equations with time-dependent linear damping, generalizing  
previous work on conformal symplectic methods, which only allows constant coefficients on the damping terms.

The aim of the present article is to 
develop these techniques for PDEs, which have damping coefficients that are necessarily functions of time 
and may be subject to external or parametric forcing.
Similar to \cite{Moore09}, we consider a multi-symplectic PDE with additional (linear) forcing and damping terms
\begin{equation}
{\bf K}z_t +{\bf L}z_x = \nabla_{z} S(z,t) -a(t){\bf K}z +{\bf F}(x,t),
\label{eq:fdms}
\end{equation}
where ${\bf F}$ is an external driving force and subscripts denote usual derivatives,
but in this case we allow the damping coefficient $a(t)$ and the smooth function $S$ to depend on time.
The corresponding variational equation takes the form
\[ {\bf K} dz_{t} +{\bf L}dz_{x} = S_{zz}(z,t) dz -a(t){\bf K}dz \]
where $S_{\zeta\zeta}$ is the Hessian.  
Taking the wedge product with $dz$ and multiplying through by 
$e^{\theta(t)}$ with $\theta(t) := \int_{0}^{t}a(s) ds$ yeilds
\begin{equation} \label{eq:mcscl}
\partial_{t} \left( e^{2\theta(t)} (dz \wedge {\bf K} dz) \right) + 
\partial_{x} \left( e^{2\theta(t)} (dz \wedge {\bf L} dz) \right) = 0,
\end{equation}
which is a multi-conformal-symplectic conservation law \cite{M16}.
We call a discretization of (\ref{eq:fdms}) {\sl multi-conformal-symplectic} 
if it satisfies a discrete version of (\ref{eq:mcscl}).
Other similar conservation laws may exist in the special case ${\bf F}(x,t) = 0$. 
Taking the inner product of the equation 
\begin{equation}
{\bf K}z_t +{\bf L}z_x = \nabla_{z} S(z,t) -a(t){\bf K}z,
\label{eq:dms}
\end{equation}
with $z_x$ \cite{MNS13}, 
and then multiplying the result through by $e^{\theta(t)}$ yields 
a conformal momentum conservation law
\begin{equation} \label{eq:momcl}
 \partial_{t} \left( e^{2\theta(t)} \left(\frac{1}{2} z^{T} {\bf K} z_x \right) \right) + 
 \partial_{x} \left( e^{2\theta(t)} \left(\frac{1}{2} z_{t}^{T} {\bf K} z +S(z,t)\right) \right) = 0.
\end{equation}
Finally, if ${\bf B}$ defines an action under which the system (\ref{eq:dms}) is invariant, 
in the special case $a(t) = 0$,
so that $({\bf B}z)^{T} \nabla S(z,t) = 0$, then solutions of the system with any $a(t)$ also 
satisfy
\begin{equation} \label{eq:lscl}
 \partial_{t} \left( e^{2\theta(t)} z^T {\bf KB} z\right) + 
 \partial_{x} \left( e^{2\theta(t)} z^T {\bf LB} z\right) = 0.
\end{equation}
We say a numerical method preserves a conservation law if it satisfies a discrete verion of it.

Overall, there are several particular PDEs with solutions that satisfy such conservation laws, 
as illustrated in Section \ref{sec:models}. 
A few numerical schemes that preserve the conformal structures of PDEs with parametric and/or external forcing and 
 damping coefficients that depend on time are presented in Section \ref{sec:SPI}.
Then, we apply specific methods to damped-driven NLS and damped Camassa-Holm equations in Section \ref{sec:NLS}, where
a multi-conformal-symplectic method is used for an equation with a norm conservation law,
experiments demonstrate advantageous numerical results.

%%%%%%%%%%%%%%%%%%%%%%%%%%%%%%%%%%%%%%%%%%%%%%%%%%%%%%%%%%%%%%%%%%%%%%%%%%%%%%%%%%%%%%%%%%%
\section{Model Problems} \label{sec:models}
In order to motivate the study of equations with the general forms (\ref{eq:fdms}) and (\ref{eq:dms}), 
we consider some examples of PDEs that may be stated in those forms
along with some of the conservation laws that are satisfied by solutions of the equations. 
Indeed, special cases (constant coefficients) of the examples presented here were also stated 
in \cite{BM17} along with their global conformal invariants, which depend on special boundary conditions.  
Here, we show that equations with variable coefficients satisfy local conservation laws that are independent
of the boundary conditions.

%%%%%%%%%%%%%%%%%%%%%%%%%%%%%%%%%%%%%%%%%%%%%%%%%%%%%%%%%%%%%%%%%%%%%%%%%%%%%%%%%%%%%%%%%%%%%%%%%%%%%%%%%%%%%%%%%%%%
\subsection{Damped wave equations}
Consider a semi-linear wave equation with time-dependent damping (cf.~\cite{W17}) given by
\[ u_{tt} = u_{xx} -2a(t)u_t -f^{\prime}(u). \]
Equations of this type are used to model propagation of waves inside matter, 
and they are the subject of many physical and mathematical studies 
(see for example \cite{DA15,DALR13,KS01,N11,NZ09,W04}). 
This equation can be rewritten as the system of equations
% \begin{equation}\label{sle}
% \begin{array}{rcl}
% -v_t -w_x &=& f^{\prime}(u) +av \\ %label{weq:sys1}\\
% u_t -p_x &=& v \\
% -p_t +u_x &=& -w +ap \\
% w_t +v_x &=& -g^{\prime}(p) %\label{weq:sys4}
% \end{array}
% \end{equation}
% which is of the form 
(\ref{eq:dms}) with $z = [u,v,w,p]^T$,
\begin{equation} \label{eq:waveJKL}
{\bf K} = \left[ \begin{array}{cc} -{\bf J} & {\bf 0} \\ {\bf 0} & -{\bf J} \end{array} \right], \qquad
{\bf L} = \left[ \begin{array}{cc} {\bf 0} & -{\bf I} \\ {\bf I} & {\bf 0} \end{array} \right], \qquad
{\bf J} = \left[ \begin{array}{cc} 0 & 1 \\ -1 & 0 \end{array} \right]
\end{equation}
and 
\[ S(z,t) =  a(t)(uv + wp) + \frac{1}{2}\left(v^2 -w^2 \right) + f(u) +a^{\prime}(t)p^2. \]
Thus, the equation satisfies conservation laws of the form (\ref{eq:mcscl}) and (\ref{eq:momcl}). 
% special constructions also have conservation laws of the form ({\ref{eq:lscl}).  
% For the case of the Klein-Gordon equation, such that
% $f^{\prime}(u) = cu$, it is a simple calculation to show that
% % \[ {\bf B} = \left[ \begin{array}{cccc} 0 & 0 & -1 & 0 \\ 0 & 0 & 0 & c \\
% % -c & 0 & 0 & 0 \\ 0 & 1 & 0 & 0 \end{array} \right] \]
% % and 
% taking the inner product of the system with ${\bf B}z :=[-w\,\,cp\,\,-\!\!cu\,\,v]^{T}$ yields the conservation law
% \begin{equation} \label{weq:ccl}
% \partial_{t} \left( e^{2\theta(t)} (vw + cup) \right)+
% \partial_{x} \left( \frac{e^{2\theta(t)}}{2} \left(v^{2} +w^{2} - c(u^2 +p^2)\right) \right)= 0.
% \end{equation}

%%%%%%%%%%%%%%%%%%%%%%%%%%%%%%%%%%%%%%%%%%%%%%%%%%%%%%%%%%%%%%%%%%%%%%%%%%%%%%%%%%%%%%%%%%%%%%%%%%%%%%%%%%%%%%%%%%%
\subsection{Generalized KdV equations}
A generalized KdV equation with time-dependent damping and dispersion, given by 
\[ u_t + u^{k-1}u_x + a(t)u + b(t)u_{xxx} = 0, \]
for $k \in \ZZ^{+}$, 
is used in studies of various physical phenomena, including coastal waves and dipole blocking 
(see for example \cite{B09,JK10,XFS06}).
The equation can be stated in the form (\ref{eq:dms}) with $z = [\phi,u,v,w]^T$
\begin{equation} \label{eq:KdVJKL}
{\bf K} = \left[ \begin{array}{cc} {\bf J} & {\bf 0} \\ {\bf 0} & {\bf 0} \end{array} \right], \qquad
{\bf L} = \left[ \begin{array}{cc} {\bf 0} & {\bf J} \\ {\bf J} & {\bf 0} \end{array} \right], \qquad
{\bf J} = \left[ \begin{array}{cc} 0 & 1 \\ -1 & 0 \end{array} \right]
\end{equation}
and $S(z,t) = \frac{1}{4b(t)}v^2 -uw +\frac{2}{k(k+1)}u^{k+1}$.
Thus, the equation clearly has conservation laws corresponding to (\ref{eq:mcscl}) and (\ref{eq:momcl}), 
and it also satisfies the conservation law 
\[ \partial_{t} \left( e^{\theta(t)} u\right) 
+\partial_{x} \left( e^{\theta(t)} \left( \frac{1}{k} u^{k} + b(t) u_{xx} \right) \right) = 0. \]

%%%%%%%%%%%%%%%%%%%%%%%%%%%%%%%%%%%%%%%%%%%%%%%%%%%%%%%%%%%%%%%%%%%%%%%%%%%%%%%%%%%%%%%%%%%%%%%%%%%%%%%%%%%%%%%%%%%%%%%
\subsection{Damped-driven NLS equations}
Several variations of the damped-driven nonlinear Schr\"{o}dinger 
(DDNLS) equation with variable coefficients, given by
\[ i\psi_{t} +\rho(t)\psi_{xx} +i\gamma(t)\psi +\sigma(t)V^{\prime}(|\psi|)\psi = 0, \]
have applications in mathematical physics (see for example \cite{KPH,SH,SMB}).
These equations can also be stated in the form (\ref{eq:dms}).
As an example, consider a damped-driven NLS equation \cite{HDY17,ChiPhys02} given by
\begin{equation} \label{eq:nls5}
 i\psi_{t} +\psi_{xx} +i\gamma\psi +ce^{i\omega t} \psi +V^{\prime}(|\psi|)\psi = 0. 
\end{equation}
Let $i\gamma\psi +ce^{i\omega t} \psi= (\alpha + i\beta)\psi$, 
where $\alpha = \alpha(t) =c\cos(\omega t)$ and  $\beta= \beta(t) =\gamma + c\sin (\omega t)$.  
Setting $\psi = p +iq$, equation (\ref{eq:nls5}) can be rewritten as
\begin{align} \label{eq:nls6}
q_t -p_{xx} &= V^{\prime}(p^2 +q^2)p +\alpha p -\beta q \\
\label{eq:nls7}
-p_t -q_{xx} &= V^{\prime}(p^2 +q^2)q +\alpha q +\beta p,
\end{align}
or, to state it in the form of equation (\ref{eq:fdms}), set $z = [p,q,v,w]^{T}$, $a(t) = \beta(t)$,
\begin{equation} \label{eq:nlsJKL}
{\bf K} = \left[ \begin{array}{cc} {\bf J} & {\bf 0} \\ {\bf 0} & {\bf 0} \end{array} \right], \qquad
{\bf L} = \left[ \begin{array}{cc} {\bf 0} & -{\bf I} \\ {\bf I} & {\bf 0} \end{array} \right], \qquad
{\bf J} = \left[ \begin{array}{cc} 0 & 1 \\ -1 & 0 \end{array} \right]
\end{equation}
and
\[ S(z,t) = \frac{1}{2} (v^2 + w^2 +V(p^2 +q^2) + \alpha(t)(p^2 +q^2) ). \]
Given this general form, it is clear that the equation has time-dependent damping/forcing terms
and satisfies a multi-conformal-symplectic conservation law.
Indeed, it also satisfies a conformal norm conservation law 
\begin{equation} \label{eq:nls9}
\partial_{t}N+\partial_{x}M=-2\beta N \quad \iff \quad
\partial_{t}(e^{2\theta(t)}N)+\partial_{x}(e^{2\theta(t)}M)=0,
\end{equation}
which is a special form of (\ref{eq:lscl}),
for $N=p^2+q^2$, $M= 2(pq_x-qp_x)$, and 
$\theta(t) = \int_{0}^{t} \beta(s)ds = \gamma t+\frac{c}{\omega}(1-\cos(\omega t)$.

Another damped NLS equation with parametric forcing \cite{BZ04,FQZ05,SB02} is given by
\begin{equation} \label{eq:DDNLSconj}
i\psi_{t} +\psi_{xx} +i\gamma\psi +ce^{i\omega t} \psi^{*} +V^{\prime}(|\psi|)\psi = 0 
\end{equation}
where $\psi^{*}$ denotes the complex conjugate of $\psi$.
Settting $\psi = p +iq$, $\alpha=c\cos(\omega t)$ and  $\beta= c\sin (\omega t)$ implies
\begin{align*}
 p_t +q_{xx} +\gamma p +\beta p -\alpha q +V^{\prime}(p^2 +q^2)q &= 0\\
 -q_t +p_{xx} -\gamma q +\beta q +\alpha p +V^{\prime}(p^2 +q^2)p &= 0.
\end{align*}
Equation (\ref{eq:DDNLSconj}) may also be cast in the form (\ref{eq:fdms}) using the matrices defined in (\ref{eq:nlsJKL})
with 
\[ S(z,t) = \frac{1}{2} (v^2 + w^2 + V(p^2 +q^2) + \alpha(t)(p^2 +q^2) + 2\beta(t) qp) \]
and $z = [p,q,v,w]^{T}$.  In this case, we have set $a(t) = \gamma$.
This system of equations does not satisfy a so called conformal norm conservation law, but it does 
satisfy a conformal momentum conservation law.
It is a straightforward calculation to show that, in this case, conservation law (\ref{eq:momcl}) is just 
\begin{equation} \label{eq:dmcl}
\partial_{t}\left( e^{2\gamma t}(pq_x -qp_x)\right) 
+ \partial_{x} \left(e^{2\gamma t} \left( p_{x}^{2} +q_{x}^{2} -pq_t +qp_t 
-\alpha(q^2 -p^2) +2\beta qp +V(p^2 +q^2) \right) \right) = 0. 
\end{equation}
% which is equivalent to 
% \[ \partial_{t}(pq_x -qp_x) 
% +\partial_{x}\left( p_{x}^{2} +q_{x}^{2} -pq_t +qp_t -\alpha(q^2 -p^2) +2\beta qp +V(p^2 +q^2) \right)
%  = -2\gamma (pq_x -qp_x). \]
%To preserve (\ref{eq:dmcl}) numerically, we can use the exponential discrete gradient method (\ref{eq:expDG}).

% A damped NLS equation with external forcing \cite{BZ04}, given by
% \[ i\psi_{t} +\psi_{xx} +i\gamma\psi +V^{\prime}(|\psi|)\psi = F(x,t), \]
% has been important in the study of deep water waves \cite{S05}.
% Conformal symplectic methods for this equation have been developed in \cite{B16}, 
% and in \cite{MNS13} for the case $F =0$ and constant $\gamma$.

%A damped discrete NLS equation with parametric forcing \cite{FQZ05} is given by 
%\[ i\dot{\psi}_{k} +\alpha \sum_{i} (\psi_{i} -\psi_{k}) +i\gamma\psi_{k} +ce^{iwt} \psi_{k}^{*} +V^{\prime}(|\psi_{k}|)\psi_{k} = 0 \]

%%%%%%%%%%%%%%%%%%%%%%%%%%%%%%%%%%%%%%%%%%%%%%%%%%%%%%%%%%%%%%%%%%%%%%%%%%%%%%%%%%%%%%%%%%%%%%%%%%%%%%%%%%%%%%%%%%%%%%%%%%%%%%%%
\subsection{Damped-driven Camassa-Holm Equation}
A dissipative Camassa-Holm equation with external forcing
\begin{align} \label{eq:dch}
u_t - u_{xxt} + 3uu_x +\gamma(t) (u - u_{xx}) - 2u_x u_{xx} -u u_{xxx} = f(x,t), \quad  u(x,0) = u_0,
\end{align}
is used to model dispersive shallow water waves on a flat bottom \cite{Shen}. 
Here, $u$ is the fluid velocity in the $x$ direction or the height of water's free surface above a flat bottom. 
%Camassa and Holm derived the conservative counterpart ($\gamma =0$) of the above equation \cite{CH93, CH94} as an asymptotic expansion in the Hamiltonian for Euler's equations. 
The conservative Camassa-Holm equation ($\gamma = 0, f =0$) is completely integrable, bi-Hamiltonian and possesses infinitely many conservation laws. 
%It has traveling wave solutions that are solitons and permanent solutions. The solitons are also called \emph{peakons} because they have peaks where the first derivatives are discontinuous. It also models wave breaking for $u_0 \in \IH^1(\IR)$ whereas KdV, another well known shallow water wave equation whose solutions are global, does not. Moreover, authors in \cite{AC00} show, among other things, the association between the conservative Camassa-Holm equation and the geodesic flow on the infinite dimensional Hilbert manifold of diffeomorphisms of the line. Conservative Camassa-Holm equation was also obtained by the authors in \cite{FF81,Fu96} as a bi-Hamiltonian generalization of the KdV equation. In fact, authors show in \cite{DG03} that conservative Camassa-Holm equation is asymptotically equivalent to KdV5 (the fifth-order integrable equation in KdV hierarchy).
%A lot of research has been done to solve the conservative Camassa-Holm equation numerically. For example, in \cite{KL05}, authors discretize conservative Camassa-Holm equation on a circle with a four stage Runge-Kutta method after taking the Fourier transform of the equation. They simulate the traveling wave solutions of the equation and peakon and cuspon interactions with the method thus obtained and observe that cuspons interact elastically. 
As a result, several authors have constructed and analyzed structure-preserving algorithms for solving the equation, including
an energy preserving numerical scheme \cite{HR06}, a multi-sympletic numerical scheme constructed through a variational approach \cite{KS00},
and multi-symplectic methods based on Euler box scheme \cite{CO08},
one of which is designed to handle conservative continuation of peakon-antipeakon collision.

Since the physical system generally exhibits some weak dissipation \cite{Shen}, 
we are interested in the dissipative ($\gamma \neq 0$) Camassa-Holm equation. 
It was shown in \cite{WY06,WY07,WY09} that the global solutions of (\ref{eq:dch}) with $f=0$
exist and decay to zero as time goes to infinity for certain initial profiles and under certain restrictions 
on the associated potential $(1-\partial_x^2)u_0$. 
Hence (\ref{eq:dch}) does not possess any traveling wave solutions although the equation does manifest 
the wave breaking phenomenon like its conservative counterpart. 
It was also shown there that the onset of blow-up is affected by the dissipative parameter $\gamma$ 
but the blow-up rate is not affected by $\gamma$. We point out the following result from these articles
$$\|u(x,t)\|_{\IH^1}^2 = e^{-2\gamma t} \|u_0(x)\|_{\IH^1}^2 \text{ for all } t\in[0,T]$$
and $u_0 \in \IH^s$, $s > \frac{3}{2}$ where $T >0$ is the maximal existence time of the solution $u$. 
In \cite{SS06}, authors establish a global well-posedness result and existence of global attractors 
in the periodic case for a closely related equation called viscous Camassa-Holm equation.

Equation (\ref{eq:dch}) can be written in conformal multi-symplectic form (\ref{eq:fdms}) with $a(t) = \gamma(t)$,
\begin{align*}
{\bf K} = \begin{bmatrix}
0 & \half & 0 & 0 & -\half \\
-\half & 0 & 0 & 0 & 0 \\
0 & 0 & 0 & 0 & 0 \\
0 & 0 & 0 & 0 & 0 \\
\half & 0 & 0 & 0 & 0 \\
\end{bmatrix}, \quad {\bf L} = \begin{bmatrix}
0 & 0 & 0 & -1 & 0 \\
0 & 0 & 1 & 0 & 0 \\
0 & -1 & 0 & 0 & 0 \\
1 & 0 & 0 & 0 & 0 \\
0 & 0 & 0 & 0 & 0 \\
\end{bmatrix}, \quad z = \begin{bmatrix}
u \\ v \\ w \\ q \\p \end{bmatrix},
\quad {\bf F}(x,t) = \begin{bmatrix}
f(x,t) \\ 0 \\ 0 \\ 0 \\ 0 \end{bmatrix}
\end{align*}
and $$S(z,t) = -\frac{u^3}{2} - \frac{up^2}{2} - uw + qp +\frac{\gamma(t)}{2}up.$$
Thus, the solutions of (\ref{eq:dch}) satisfy a multi-conformal-symplectic conservation law, 
and it is a straightforward calculation to show that the equation satisfies other conservation laws, given by 
\[ \partial_{t} \left( e^{\theta(t)} u \right) + \partial_{x} \left( e^{\theta(t)} \left(  
\frac{3}{2}u^2 -\frac{1}{2} u_{x}^{2} -u_{xt} -uu_{xx} -\gamma u_x \right) \right) = 0 \] 
and 
\[ \partial_{t} \left( \frac{e^{\theta(t)}}2 \left( u^2 + u_{x}^{2} \right) \right) 
+ \partial_{x} \left( e^{\theta(t)} \left( u^3 -u^2 u_{xx} -uu_{xt} -uu_x \right) \right) =0 \]
in the absence of external forcing.  
% , equation (\ref{eq:dch}) satisfies the following conformal properties
% % \begin{align} \label{eq:dch-cp}
% % \int u\,dx, \quad \int (u^2 +u_x^2)\,dx,
% % \end{align}
% % i.e.
% \begin{align*}
% \frac{d}{dt} \int u\,dx = -\gamma \int u\,dx, \quad \frac{d}{dt} \int (u^2 +u_x^2)\,dx = -2\gamma \int (u^2 +u_x^2)\,dx.
% \end{align*}
% and 

%%%%%%%%%%%%%%%%%%%%%%%%%%%%%%%%%%%%%%%%%%%%%%%%%%%%%%%%%%%%%%%%%%%%%%%%
\section{Structure-Preserving Integrators} \label{sec:SPI}
%%%%%%%%%%%%%%%%%%%%%%%%%%%%%%%%%%%%%%%%%%%%%%%%%%%%%%%%%%%%%%%%%
Given that there are many interesting PDEs that may be formulated in the general forms 
(\ref{eq:fdms}) or (\ref{eq:dms}), it is desirable to construct numerical integrators that 
preserve one or more of the conservation laws that are satisfied by solutions of the equations.
There are two main approaches for constructing conformal symplectic methods for solving conformal Hamiltonian ODEs.  
The first, proposed by McLachlan and Quispel \cite{McLQ01,McLQ02}, 
is a splitting approach in which the vector field is split into a Hamiltonian part and 
a linearly damped part.  Then the exact flow map of the damped equation is composed with a symplectic flow map that 
is used to solve the conservative equation, resulting in a flow map that is conformal symplectic.  
Another approach relies on choosing the coefficient functions of an exponential Runge-Kutta method 
(or a partitioned exponential Runge-Kutta method) 
in a way a guarantees preservation of conformal symplecticity \cite{BM17}.
Each approach has advantages, and both lead to some of the same methods that have been studied in previous works.
From another point of view, a Birkhoffian formulation of the PDE  
also produces some of the same methods \cite{SQWS10}.

Each of these approaches for constructing conformal symplectic integrators 
is based on an underling symplectic method, 
but in each case the resulting conformal symplectic method is not guaranteed 
to maintain the order of the underlying method (cf. \cite{MNS13,SS05}).
Though splitting methods can be constructed in such a way that the order is known, 
a simpler approach is to apply a change of variables, similar to a Lawson transformation \cite{Law67}, 
and use multi-symplectic method on the transformed equation.  
This approach was used to construct structure-preserving integrating-factor methods for 
systems of ODEs \cite{BM17}.  For the PDE (\ref{eq:fdms}) the approach can be delineated as follows.
Introduce the change of variables $\zeta = e^{\theta(t)} z$ 
where $\theta(t) := \int_{0}^{t}a(s) ds$, and write equation (\ref{eq:fdms}) as
\begin{equation} \label{eq:expMS}
 {\bf K} \zeta_{t} +{\bf L}\zeta_{x} =\nabla \tilde{S}(\zeta,t) +{\bf G}(x,t), 
\end{equation}
where ${\bf G}(x,t) =e^{\theta(t)} {\bf F}(x,t)$ 
and $\tilde{S}(\zeta,t) = e^{\theta(t)} S(e^{-\theta(t)} \zeta,t)$.
Then, apply any number of structure-preserving algorithms to this transformed 
equation, which is the so-called underlying method. 
It is expected that transforming back to the original variables 
will result in methods that preserve conformal versions of the properties, provided they exist.
With this approach the order of the underlying method is maintained in general.

Here, we present discretization methods that preserve one or more of the conservation laws 
(\ref{eq:mcscl}), (\ref{eq:momcl}), and (\ref{eq:lscl}).  
Due to the structure of the equations and the damping/driving terms, it is 
essential to use an exponential method for the time discretization, if structure preservation is important. 
Indeed, there many advantageous possibilities provided by a special class of exponential Runge-Kutta methods \cite{BM17}.
For simplicity, a second-order exponential method \cite{BFM15}, 
which reduces to the implicit midpoint rule in the case $a=0$ and is denoted by a Butcher-like tableau
\begin{equation}
    \label{fig:ERK-Lawson}
    \begin{array}{l | c | l}
      \rule{0pt}{2,3ex} \frac{1}{2}      \quad & \frac{1}{2} & e^{-\int_{t_n}^{t_{n+1/2}} a(s) ds} \\
\hline
      \rule{0pt}{3,3ex}        \quad & e^{-\int_{t_{n+1/2}}^{t_{n+1}} a(s) ds} & e^{-\int_{t_n}^{t_{n+1}} a(s) ds}
    \end{array}
\end{equation}
is used to discretize the PDE in time in the following subsections.  

Many standard spatial discretizations can be used.  In fact,
we expect that the spatial discretization will guarantee preservation of the desired conservation law, as long 
as it preserves the same conservation law in the absence of damping/driving forces, meaning  
finite-element, pseudospectral, and finite-difference, methods are all possibilities. 
As a introduction to structure-preserving discetizations for PDEs with time-dependent damping/driving forces, we 
present only three of the most basic finite-difference spatial discretizations in the following subsections.
These discretizations are similar to the methods proposed in \cite{MNS13}, but with two important differences.  
First, the methods of \cite{MNS13} are first-order in time for general nonlinear equations, but 
the methods proposed here are second-order.  Second, the methods of \cite{MNS13} assume that $a(t)$ is 
constant, and that $S$ does not explicitly depend on time.

To begin, let $z_{n}^{i}$ denote and approximation of $z(x_n,t_i)$, where 
$x_{n+1} = x_n +\dx$ and $t_{i+1} = t_i +\dt$.  Also, define
\[ z_{n+1/2}^{i} = \frac{z_{n+1}^{i} +z_{n}^{i}}{2}, \qquad
\delta_{x}^{-} \zni = \frac{\zni -z_{n-1}^{i}}{\dx}, \qquad \textup{and} \qquad 
 \delta_{x}^{+} \zni = \frac{z_{n+1}^{i} -\zni}{\dx}, \]
along with 
\[ A_{t_i}^{a} \zni= \frac{1}{2} \left( e^{\int_{t_{i+1/2}}^{t_{i+1}} a(s)ds}  z_{n}^{i+1} 
+  e^{-\int_{t_{i}}^{t_{i+1/2}} a(s)ds}  z_{n}^{i} \right)\]
and 
\[ D_{t_i}^{a} \zni= \frac{1}{\dt} \left( e^{\int_{t_{i+1/2}}^{t_{i+1}} a(s)ds}  z_{n}^{i+1} 
-  e^{-\int_{t_{i}}^{t_{i+1/2}} a(s)ds}  z_{n}^{i} \right).  \]
To show that discretizations using these operators are structure-preserving, it is useful 
to have a discrete product rule
\begin{equation} \label{eq:prodrule}
D_{t_i}^{2a} \left[ (z^i)^{T} y^i \right] = 
\left( D_{t_i}^{a} z^i \right)^{T} A_{t_i}^{a} y^i + \left( A_{t_i}^{a} z^i \right)^{T} D_{t_i}^{a} y^i
\end{equation}
which follows from
\begin{align*}
D_{t_i}^{2a} \left[ (z^i)^{T} y^i \right] =&
\frac{1}{\dt}\left[  e^{\int_{t_{i+1/2}}^{t_{i+1}} 2a(s)ds} (z^{i+1})^{T} y^{i+1}  
-e^{-\int_{t_{i}}^{t_{i+1/2}} 2a(s)ds}  (z^i)^{T} y^i \right] \\
=& \frac{1}{2\dt}\left( e^{\int_{t_{i+1/2}}^{t_{i+1}} a(s)ds} z^{i+1} 
- e^{-\int_{t_{i}}^{t_{i+1/2}} a(s)ds} z^i \right)^{T}  \left( e^{\int_{t_{i+1/2}}^{t_{i+1}} a(s)ds} y^{i+1} 
+ e^{-\int_{t_{i}}^{t_{i+1/2}} a(s)ds} y^i \right) \\
& + \frac{1}{2\dt}\left( e^{\int_{t_{i+1/2}}^{t_{i+1}} a(s)ds} z^{i+1} 
+ e^{-\int_{t_{i}}^{t_{i+1/2}} a(s)ds} z^i \right)^{T}  \left( e^{\int_{t_{i+1/2}}^{t_{i+1}} a(s)ds} y^{i+1} 
- e^{-\int_{t_{i}}^{t_{i+1/2}} a(s)ds} y^i \right). 
\end{align*}
In addition, notice that
\begin{align*}
 \delta_{t}^{+} \left( e^{\int_{0}^{t_i} a(s)ds} y^i \right) =& 
 \frac{1}{\dt} \left( e^{\int_{0}^{t_{i+1}} a(s)ds} y^{i+1} - e^{\int_{0}^{t_i} a(s)ds} y^i \right)\\
 =& \frac{1}{\dt} e^{\int_{0}^{t_{i+1/2}} a(s)ds} \left( 
  e^{\int_{t_{i+1/2}}^{t_{i+1}} a(s)ds} y^{i+1} - e^{-\int_{t_i}^{t_{i+1/2}} a(s)ds} y^i \right). 
\end{align*}
which leads to a useful relationship between $D_{t_i}^{a}$ 
and the standard forward difference operator given by 
\begin{equation} \label{eq:expdiff}
 \delta_{t}^{+} \left( e^{\theta_{i}} y^i \right) = 
 e^{\theta_{i+1/2}} D_{t_i}^{a} y^i 
\end{equation}
where we have intruduced $\theta_{i} := \int_{0}^{t_i} a(s)ds$ for the sake of convenience.  

%%%%%%%%%%%%%%%%%%%%%%%%%%%%%%%%%%%%%%%%%%%%%%%%%%%%%%%%%%%%%%%%%%%%%%%%%%%%%%%%%%%%%%%%%%%%%%%%%%%%%%%%%%%%%%%%%%%%
\subsection{Symplectic Euler spatial discretization}
Let us rewrite \eqref{eq:fdms} as follows
\begin{align} \label{eq:fdms-1}
{\bf K}z_t +{\bf L_+}z_x +{\bf L_-}z_x = \nabla_{z} S(z,t) -a(t){\bf K}z +{\bf F}(x,t),
\end{align}
where ${\bf L_+}+{\bf L_-}={\bf L}$ and ${\bf L_+} = -{\bf L}_-^T$. 
% Then discretizing (\ref{eq:fdms-1}) in space we get
% \begin{align*}
% {\bf K}\partial_t z_n +{\bf L_+}\delta_x^+ z_n +{\bf L_-}\delta_x^+ z_{n-1} = \nabla_{z} S(z_n,t) -a(t){\bf K}z_n +{\bf F}(x_n,t),
% \end{align*}
% where $z_n$ denotes the semi-discretized solution vector $\{z_n^i\}_{n}$. We now discretize this semi-discretized equation in time to get the fully discrete equation
Now, consider the discretization method
\begin{align} \label{eq:embs}
{\bf K}D_{t_i}^{a} z_n^i +{\bf L_+}A_{t_i}^{a} \delta_x^+ z_n^i +{\bf L_-}A_{t_i}^{a} \delta_x^- z_{n}^i 
= \nabla_{z} S(A_{t_i}^{a}z_n^i,t_{i+1/2}) +{\bf F}(x_n,t_{i+1/2}).
\end{align}
Here, the \emph{underlying method} (i.e.~the method that results when $a=0$ and $F=0$) 
is a mixed box scheme with implicit midpoint in time and symplectic Euler in space (cf. \cite[Chapter 4]{thesis}).

We now show that (\ref{eq:embs}) is a multi-conformal-symplectic scheme. The variational equation associated with this scheme is
\begin{align*}
{\bf K}D_{t_i}^{a} dz_n^i +{\bf L_+}A_{t_i}^{a} \delta_x^+ dz_n^i +{\bf L_-}A_{t_i}^{a} \delta_x^+ dz_{n-1}
= S_{zz}(A_{t_i}^{a}z_n^i,t_{i+1/2})A_{t_i}^{a}dz_n^i.
\end{align*}
where $S_{zz}$ is the Hessian of $S$.
Taking the wedge product of $A_{t_i}^{a} dz_n^i$ with this equation we get
\begin{align} \label{eq:embs-vareq}
A_{t_i}^{a} dz_n^i \wedge {\bf K}D_{t_i}^{a} dz_n^i +A_{t_i}^{a} dz_n^i \wedge{\bf L_+}A_{t_i}^{a} \delta_x^+ dz_n^i 
+A_{t_i}^{a} dz_n^i \wedge{\bf L_-}A_{t_i}^{a} \delta_x^+ dz_{n-1} =0,
%\notag \\= S_{zz}(A_{t_i}^{a}z_n^i,A_{t_i}^{a}t_i)(A_{t_i}^{a} dz_n^i \wedge A_{t_i}^{a}dz_n^i).
\end{align}
because $A_{t_i}^{a} dz_{n}^i \wedge S_{zz}(A_{t_i}^{a}z_n^i,t_{i+1/2}) A_{t_i}^{a} dz_{n}^i = 0$,
due to the symmetry of $S_{zz}$.
Using the discrete product rule, we get
\[ A_{t_i}^{a} dz_n^i \wedge {\bf K}D_{t_i}^{a} dz_n^i = D_{t_i}^{2a} (\half (dz_n^i \wedge {\bf K} dz_n^i)),\] 
and
\[ A_{t_i}^{a} dz_n^i \wedge{\bf L_+}A_{t_i}^{a} \delta_x^+ dz_n^i +A_{t_i}^{a} dz_n^i \wedge{\bf L_-}A_{t_i}^{a} \delta_x^+ dz_{n-1} 
= \delta_x^+ (A_{t_i}^{a} dz_{n-1}^i \wedge {\bf L}_+ A_{t_i}^{a} dz_{n}^i ).\]
These equations, along with (\ref{eq:embs-vareq}), give
\begin{align*}
D_{t_i}^{2a} (\half (dz_n^i \wedge {\bf K} dz_n^i)) + \delta_x^+ (A_{t_i}^{a} dz_{n-1}^i \wedge {\bf L}_+ A_{t_i}^{a} dz_{n}^i ) = 0,
\end{align*}
which is a discrete version of the multi-conformal-symplectic conservation law. This proves that (\ref{eq:embs}) is a multi-conformal-symplectic scheme.

The scheme does not, in general, preserve the other conservation laws, (\ref{eq:momcl}) and (\ref{eq:lscl}).  
However, the scheme does preserve such conformal conservation laws in special cases, as is demonstrated in the following section.
It is also important to notice that the external forcing plays no role in the variational equation, so 
any multi-conformal-symplectic scheme will also preserve the conservation law (\ref{eq:mcscl}) in the presence of external forcing.

%%%%%%%%%%%%%%%%%%%%%%%%%%%%%%%%%%%%%%%%%%%%%%%%%%%%%%%%%%%%%%%%%%%%%%%%%%%%%%%%%%%%%%%%%%%%%%%%%%%%%%%%%%%%%%%%%%%%
\subsection{Implicit midpoint spatial discretization}
A second-order method for solving (\ref{eq:fdms}) is given by
\begin{equation} \label{eq:expbox}
{\bf K} D_{t_i}^{a} z_{n+1/2}^{i} + {\bf L} A_{t_i}^{a} \delta_{x}^{+} \zni 
 = \nabla S\left( A_{t_i}^{a} z_{n+1/2}^{i} \right) + {\bf F}(x_{n+1/2},t_{i+1/2}), 
 \end{equation}
and the underlying method in this case is the well-known Preissmann box scheme.  
To show that (\ref{eq:expbox}) is a multi-conformal-symplectic scheme, 
we take the wedge product of the corresponding variational equation with $A_{t_n}^{a} dz_{n+1/2}^{i}$ to get
\[ A_{t_i}^{a} dz_{n+1/2}^{i} \wedge {\bf K} D_{t_i}^{a} z_{n+1/2}^{i} + 
A_{t_i}^{a} dz_{n+1/2}^{i} \wedge {\bf L} A_{t_i}^{a} \delta_{x}^{+} \zni = 0. \] 
Since the wedge product is skew-symmetric, we may use the discret product rule to write this as 
\[ D_{t_i}^{2a} \left( dz_{n+1/2}^{i} \wedge {\bf K} z_{n+1/2}^{i} \right) + 
\delta_{x}^{+} \left( A_{t_i}^{a} dz_{n}^{i} \wedge {\bf L} A_{t_i}^{a} \zni \right) = 0. \]
Then, formula (\ref{eq:expdiff}) implies
\[ \delta_{t}^{+} \left( e^{2\theta_{i}} dz_{n+1/2}^{i} \wedge {\bf K} z_{n+1/2}^{i} \right) + 
\delta_{x}^{+} \left( e^{2\theta_{i+1/2}} A_{t_i}^{a} dz_{n}^{i} \wedge {\bf L} A_{t_i}^{a} \zni \right) =0, \]
which is a discrete version of the multi-conformal-symplectic conservation law (\ref{eq:mcscl}). 

This method also satisfies a discrete version of the conservation law (\ref{eq:lscl}) when $F=0$.  
To show this, take the inner product of (\ref{eq:expbox}) with ${\bf B} A_{t_i}^{a} z_{n+1/2}^{i}$, 
and rearrange the terms to get 
\[  \left( A_{t_i}^{a} z_{n+1/2}^{i} \right)^{T} {\bf KB} D_{t_i}^{a} z_{n+1/2}^{i} +
\left( A_{t_i}^{a} z_{n+1/2}^{i} \right)^{T} {\bf LB} A_{t_i}^{a} \delta_{x}^{+} \zni = 0. \]
Then, the discrete product rule implies 
\[ D_{t_i}^{2a} \left( \left( z_{n+1/2}^{i} \right)^{T} {\bf KB} z_{n+1/2}^{i} \right) +
\delta_{x}^{+} \left( \left( A_{t_i}^{a} \zni \right)^{T} {\bf LB} A_{t_i}^{a}  \zni \right) = 0  \] 
and property (\ref{eq:expdiff}) yields the following discrete version of (\ref{eq:lscl})
\[ \delta_{t}^{+} \left( e^{2\theta_{i}} \left( z_{n+1/2}^{i} \right)^{T} {\bf KB} z_{n+1/2}^{i} \right) +
\delta_{x}^{+} \left( e^{2\theta_{i+1/2}} 
\left( A_{t_i}^{a} \zni \right)^{T} {\bf LB} A_{t_i}^{a}  \zni \right) = 0.  \] 

%%%%%%%%%%%%%%%%%%%%%%%%%%%%%%%%%%%%%%%%%%%%%%%%%%%%%%%%%%%%%%%%%%%%%%%%%%%%%%%%%%%%%%%%%%%%%%%%%%%%%%%%%%%%%% 
\subsection{Discrete gradient spatial discretizations}
In some cases it may be more desirable to preserve a conformal momentum conservation law, rather than 
the multi-conformal-symplectic conservation law.  Thus, we also propose discrete gradient discretizations. 
Define a discrete gradient that satisfies the following properties
\[ \lim_{\hat{z} \to z} \overline{\nabla}S(\hat{z},z,t) = \overline{\nabla}S(z,z,t) = \nabla S(z,t) 
 \qquad \textup{and} \qquad \overline{\nabla}S(\hat{z},z,t)^{T} (\hat{z} -z) = S(\hat{z},t) - S(z,t).
\]
A second-order discrete gradient method for solving (\ref{eq:dms}) is given by
\begin{equation} \label{eq:expDG}
{\bf K} D_{t_i}^{a} z_{n+1/2}^{i} + {\bf L} A_{t_i}^{a} \delta_{x}^{+} \zni 
 = \overline{\nabla} S\left( A_{t_i}^{a} z_{n+1}^{i}, A_{t_i}^{a} z_{n}^{i}, t_{i+1/2} \right). 
\end{equation}
To show that this method satisfies a discrete version of the momentum conservation law (\ref{eq:momcl}),
take the inner product of (\ref{eq:expDG}) with $A_{t_i}^{a} \delta_{x}^{+} \zni$ to get
\[ \left( A_{t_i}^{a} \delta_{x}^{+} \zni \right)^{T} {\bf K} D_{t_i}^{a} z_{n+1/2}^{i} 
= \frac{1}{\dx} \left( 
S\left( A_{t_i}^{a} z_{n+1}^{i}, t_{i+1/2} \right) - S\left( A_{t_i}^{a} z_{n}^{i}, t_{i+1/2} \right) \right)
= \delta_{x}^{+} S\left( A_{t_i}^{a} z_{n}^{i}, t_{i+1/2} \right). \]
Then, using the discrete product rule gives
\[ D_{t_i}^{2a} \left( \frac{1}{2} \left( \delta_{x}^{+} \zni \right)^{T} {\bf K}  z_{n+1/2}^{i} \right)
+ \delta_{x}^{+} \left( \frac{1}{2} \left( A_{t_i}^{a} \zni \right)^{T} {\bf K} D_{t_i}^{a} z_{n}^{i} \right)
= \delta_{x}^{+} S\left( A_{t_i}^{a} z_{n}^{i}, t_{i+1/2} \right).\]
Rearranging terms and using property (\ref{eq:expdiff}) yields
\[ \delta_{t}^{+} \left( e^{2\theta_{i}} 
\frac{1}{2} \left( z_{n+1/2}^{i} \right)^{T} {\bf K} \delta_{x}^{+} \zni \right)
+ \delta_{x}^{+} \left( e^{2\theta_{i+1/2}} \left( \frac{1}{2} 
\left( D_{t_i}^{a} z_{n}^{i} \right)^{T} {\bf K} A_{t_i}^{a} \zni 
+ S\left( A_{t_i}^{a} z_{n}^{i}, t_{i+1/2} \right) \right) \right) = 0.\]

%%%%%%%%%%%%%%%%%%%%%%%%%%%%%%%%%%%%%%%%%%%%%%%%%%%%%%%%%%%%%%%%%%%%%%%%
\section{Numerical Experiments} \label{sec:NLS}
%%%%%%%%%%%%%%%%%%%%%%%%%%%%%%%%%%%%%%%%%%%%%%%%%%%%%%%%%%%%%%%%%
Overall, there are many PDE models that could be characterized under the premise of this article, along with 
many competitive discretizations for each model, and it is not possible to give a thorough numerical exposition here.   
But, to illustrate the usefulness of the proposed discretizations, we perform numerical experiments 
for a damped-driven NLS equation (\ref{eq:nls5}), and a damped-driven Camassa-Holm equation (\ref{eq:dch}). 

%%%%%%%%%%%%%%%%%%%%%%%%%%%%%%%%%%%%%%%%%%%%%%%%%%%%%%%%%%%%%%%%%%%%%%%%%%%%%%%%%%%%%%%%%%%%%%%%
\subsection{Damped-driven NLS equation}
In addition to a multi-conformal-symplectic conservation law, 
one can also show that discretization (\ref{eq:embs}) (with $F=0$) preserves
the norm conservation law (\ref{eq:nls9}) for equation (\ref{eq:nls5}). 
Writing the discretization for the system (\ref{eq:nls6})-(\ref{eq:nls7}) gives 
\begin{align} 
 D_{t_i}^{\beta} q_{n}^{i} - A_{t_i}^{\beta} \delta_{x}^{2} p_{n}^{i} 
 -V^{\prime}\left( (A_{t_i}^{\beta} p_{n}^{i})^{2} + (A_{t_i}^{\beta} q_{n}^{i})^{2} \right) A_{t_i}^{\beta} p_{n}^{i} -\alpha(t_{i+1/2}) A_{t_i}^{\beta} p_{n}^{i} &= 0 
 \label{eq:nlsdisc1} \\
D_{t_i}^{\beta} p_{n}^{i} + A_{t_i}^{\beta} \delta_{x}^{2} q_{n}^{i} 
 +V^{\prime}\left( (A_{t_i}^{\beta} p_{n}^{i})^{2} + (A_{t_i}^{\beta} q_{n}^{i})^{2} \right) A_{t_i}^{\beta} q_{n}^{i} +\alpha(t_{i+1/2}) A_{t_i}^{\beta} q_{n}^{i} &= 0,
 \label{eq:nlsdisc2}
\end{align}
where $\delta_{x}^{2} := \delta_{x}^{+}\delta_{x}^{-} = \delta_{x}^{-}\delta_{x}^{+}$.
Multiplying the first equation through by $A_{t_i}^{\beta} q_{n}^{i}$ and the second equation by $A_{t_i}^{\beta} p_{n}^{i}$ and adding them together implies
\[  A_{t_i}^{\beta} q_{n}^{i} D_{t_i}^{\beta} q_{n}^{i} - A_{t_i}^{\beta} q_{n}^{i} A_{t_i}^{\beta} \delta_{x}^{2} p_{n}^{i} 
+ A_{t_i}^{\beta} p_{n}^{i} D_{t_i}^{\beta} p_{n}^{i} + A_{t_i}^{\beta} p_{n}^{i} A_{t_i}^{\beta} \delta_{x}^{2} q_{n}^{i} = 0.\] 
Then, using the discrete product rules provides
\[ D_{t_i}^{2\beta} \left( (p_{n}^{i})^{2} +(q_{n}^{i})^{2} \right) +
\frac{1}{\dx} \delta_{x}^{+} \left( A_{t_i}^{\beta}p_{n-1}^{i} A_{t_i}^{\beta} q_{n}^{i} - A_{t_i}^{\beta}p_{n}^{i} A_{t_i}^{\beta} q_{n-1}^{i} \right) =0, \]
or equivalently
\[ \delta_{t}^{+} \left( e^{2\theta_{i}} \left( (p_{n}^{i})^{2} +(q_{n}^{i})^{2} \right) \right) +
\delta_{x}^{+} \left( \frac{e^{2\theta_{i+1/2}}}{\dx} \left( A_{t_i}^{\beta}p_{n-1}^{i} A_{t_i}^{\beta} q_{n}^{i} -A_{t_i}^{\beta}p_{n}^{i} A_{t_i}^{\beta} q_{n-1}^{i} \right)\right)=0, \]
which is a discrete version of the norm conservation law (\ref{eq:nls9}) for the damped-driven nonlinear Schr\"{o}dinger equation.

\begin{figure}[htbp]
\centering
% \subfloat[]{
%   \includegraphics[width=0.5\textwidth]{exp_mixed_box_plane_wave_xy.eps}
% }
% \subfloat[]{
%   \includegraphics[width=0.5\textwidth]{exp_mixed_box_soliton_pair.eps}
% }
\hspace{0mm}
\subfloat[]{
  \includegraphics[width=0.5\textwidth]{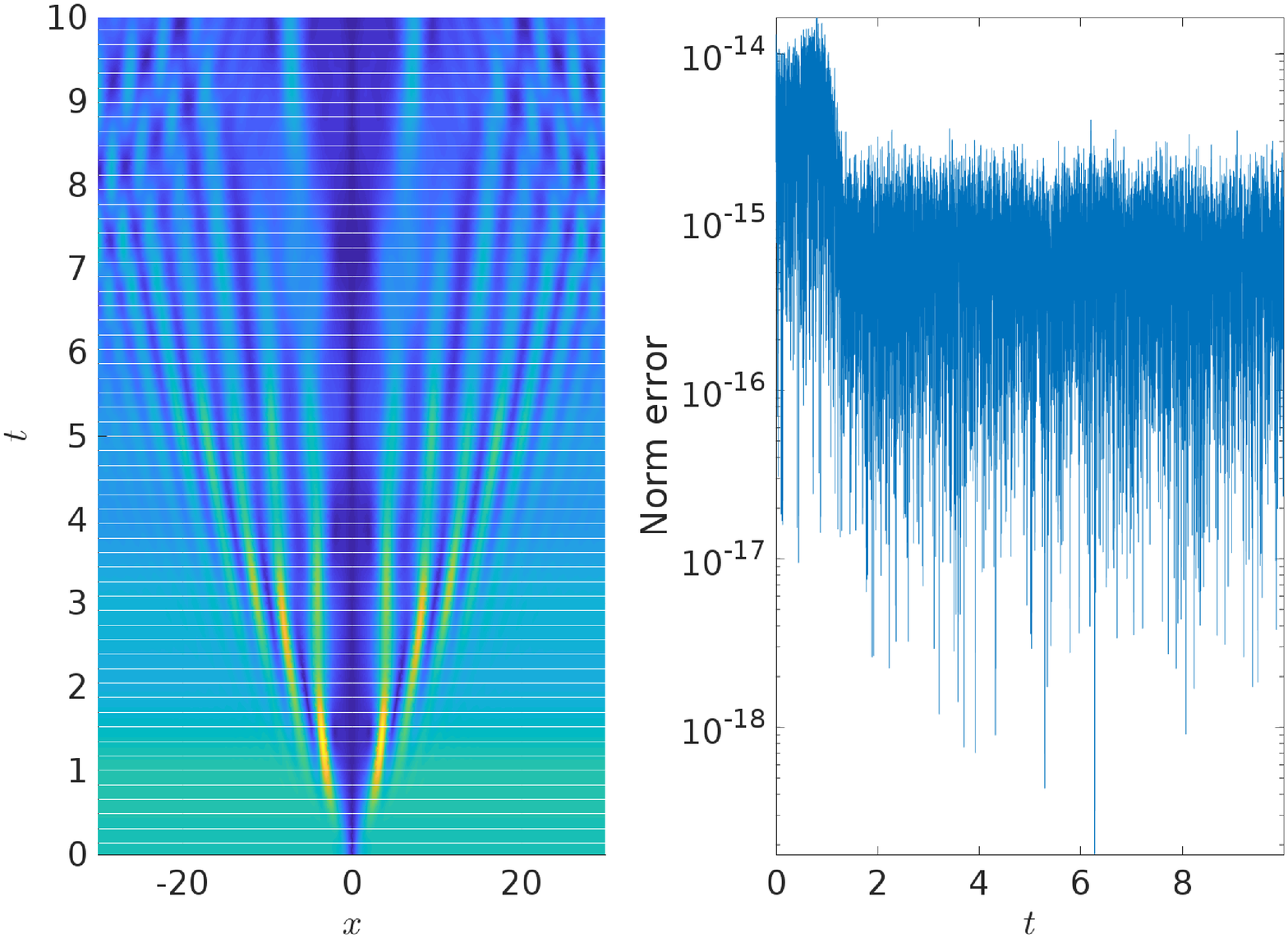}
}
\subfloat[]{
  \includegraphics[width=0.5\textwidth]{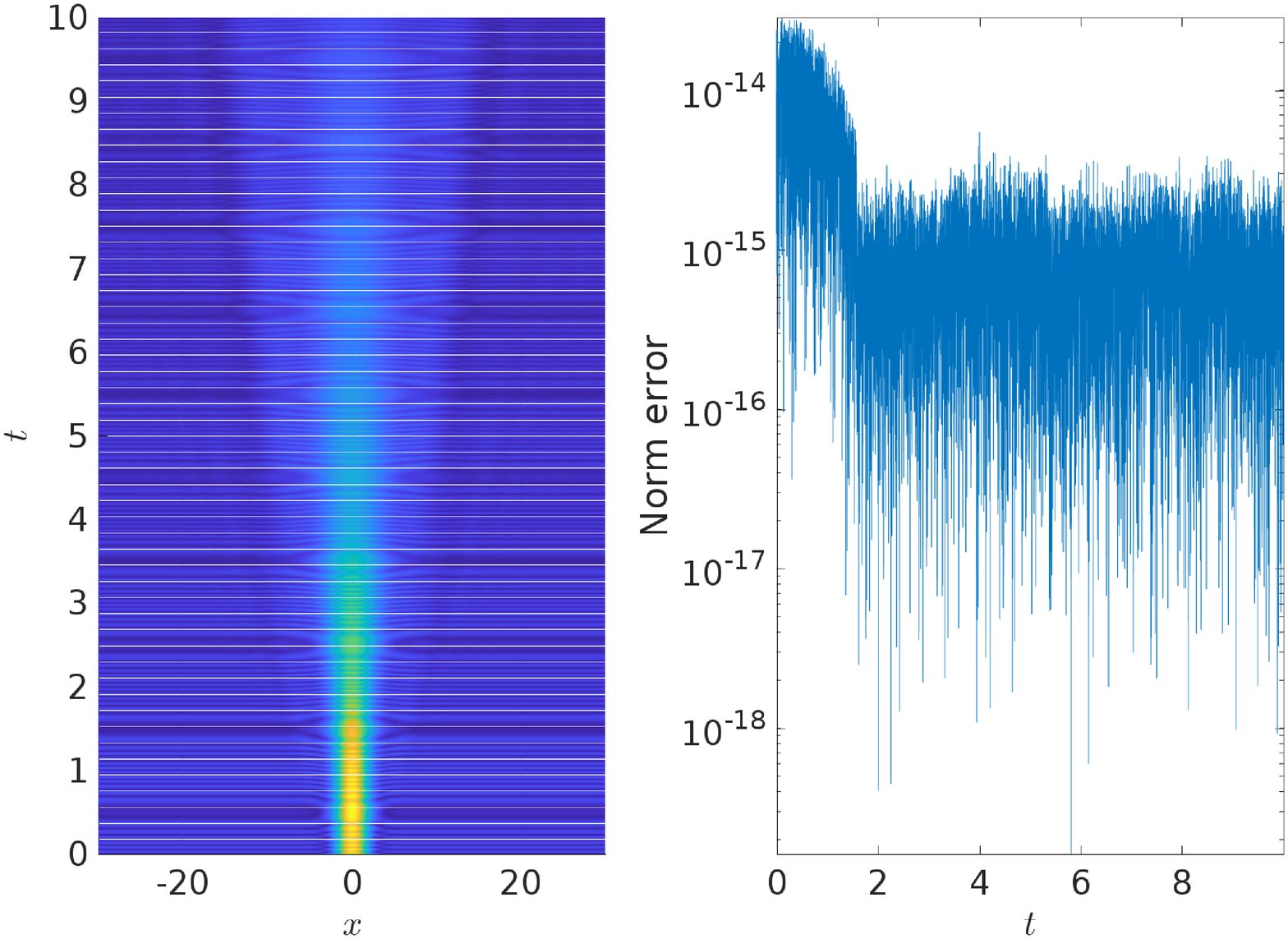}
}
\caption{Solutions (in the form $|\psi|$) of (\ref{eq:nls5}) by the method of (\ref{eq:nlsdisc1})-(\ref{eq:nlsdisc2}) next to the error in the norm.  
(a) $\psi(x,0) = \tanh(x)$ with anti-periodic boundary conditions,  
(b) $\psi(x,0) = \sqrt{\frac{2}{3\sqrt \pi}}\exp(-(2x/3)^2/2)$ with periodic boundary conditions.}
 \label{fig:embs-soln}
\end{figure}

In all of the following experiments we set
$$\dx = 0.1, ~\dt = 0.001, ~x \in [-30,30], ~t \in [0,10],~ \alpha(t) = -0.2\cos(\pi t), ~ \textup{and} ~ \beta(t) = 0.1 -0.2\sin(\pi t)$$
so that $\beta(t)$ fluctuates trigonometrically between -0.1 and 0.3. Figure \ref{fig:embs-soln} shows propagation of two initial wave profiles with suitable boundary conditions. 
One can see that the amplitude of the propagating wave increases with negative values of the parameter $\beta$ and decreases with its positive values. This causality becomes less pronounced over time as the initial profile evolves and gets dispersed over the entire spatial domain.
Propagating a dark soliton with the method of (\ref{eq:embs}) shows the solution disintegrating into a wave train \cite{KL98}(left of Figure \ref{fig:embs-soln}). 
%On the other hand, a soliton initial condition without spatial translation is seen to have stable propagation (right of Figure \ref{fig:embs-soln}). 
% We have used the following initial conditions for the respective plots. 
% All initial conditions use periodic boundary conditions except the dark soliton (c), which uses anti periodic boundary conditions.
% \begin{enumerate}
% \item $\psi(x,0) = 0.5.$
% \item $\psi(x,0) = \exp(8\ti x)\sech(x+5) +1.5\exp(-7\ti x)\sech(1.5(x-5)).$
% \item $\psi(x,0) = \tanh(x).$
% \item $\psi(x,0) = \sqrt (\frac{2}{3\sqrt(pi)})\exp(-(2/3x)^2/2).$
% \end{enumerate}
The solutions illustrate preservation of norm, by plotting the  
norm error, which is defined to be
$$\text{Norm error }= \log \left( \frac{\sum_n |\psi_n^{i+1}|^2}{\sum_n |\psi_n^{i}|^2} \right) -4(\theta_{i+1} -\theta_i).$$

\begin{figure}[htbp]
\centering
\subfloat[]{
   \includegraphics[width=0.5\textwidth]{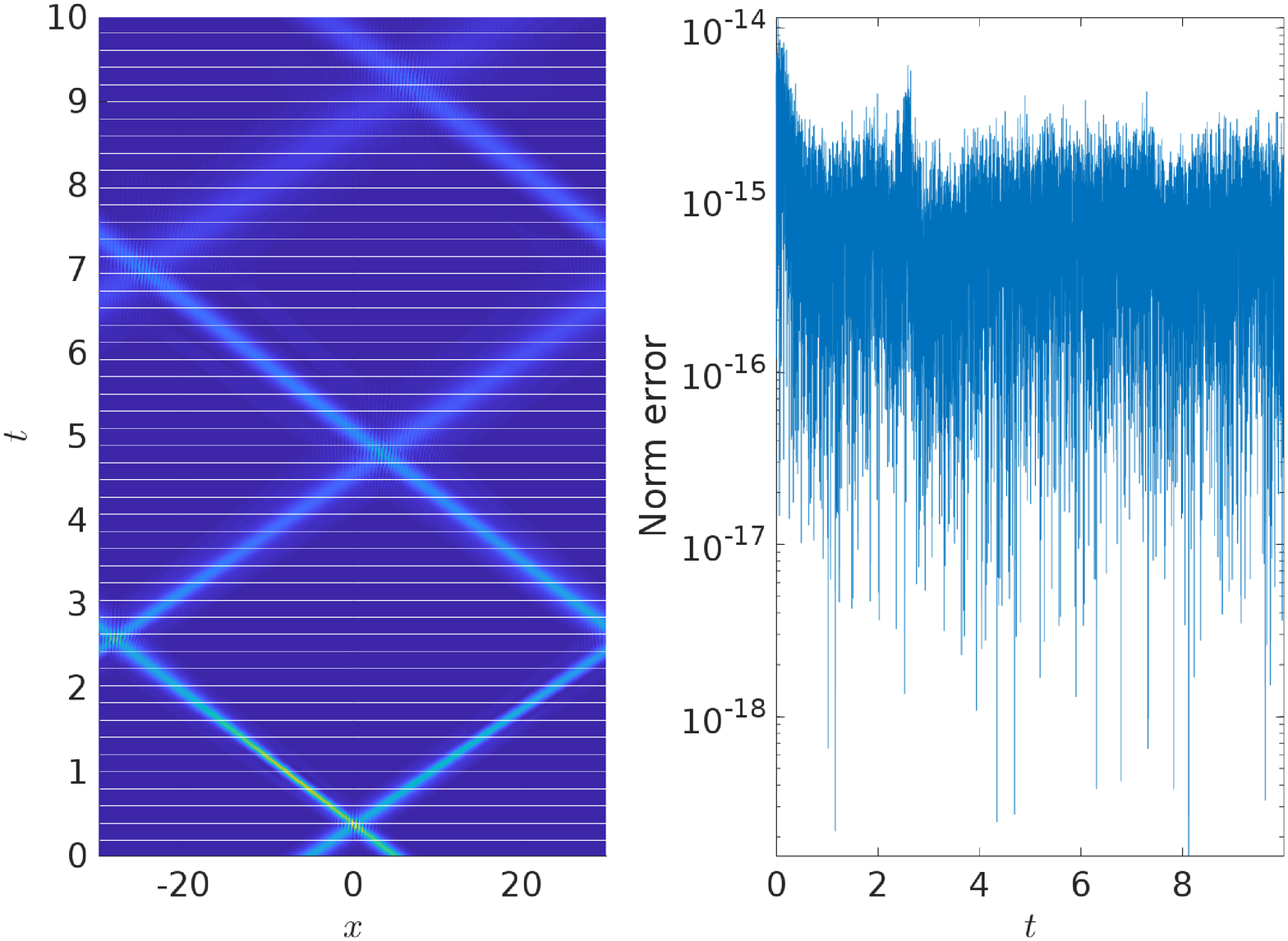}
 }
 \subfloat[]{
\includegraphics[width=0.5\textwidth]{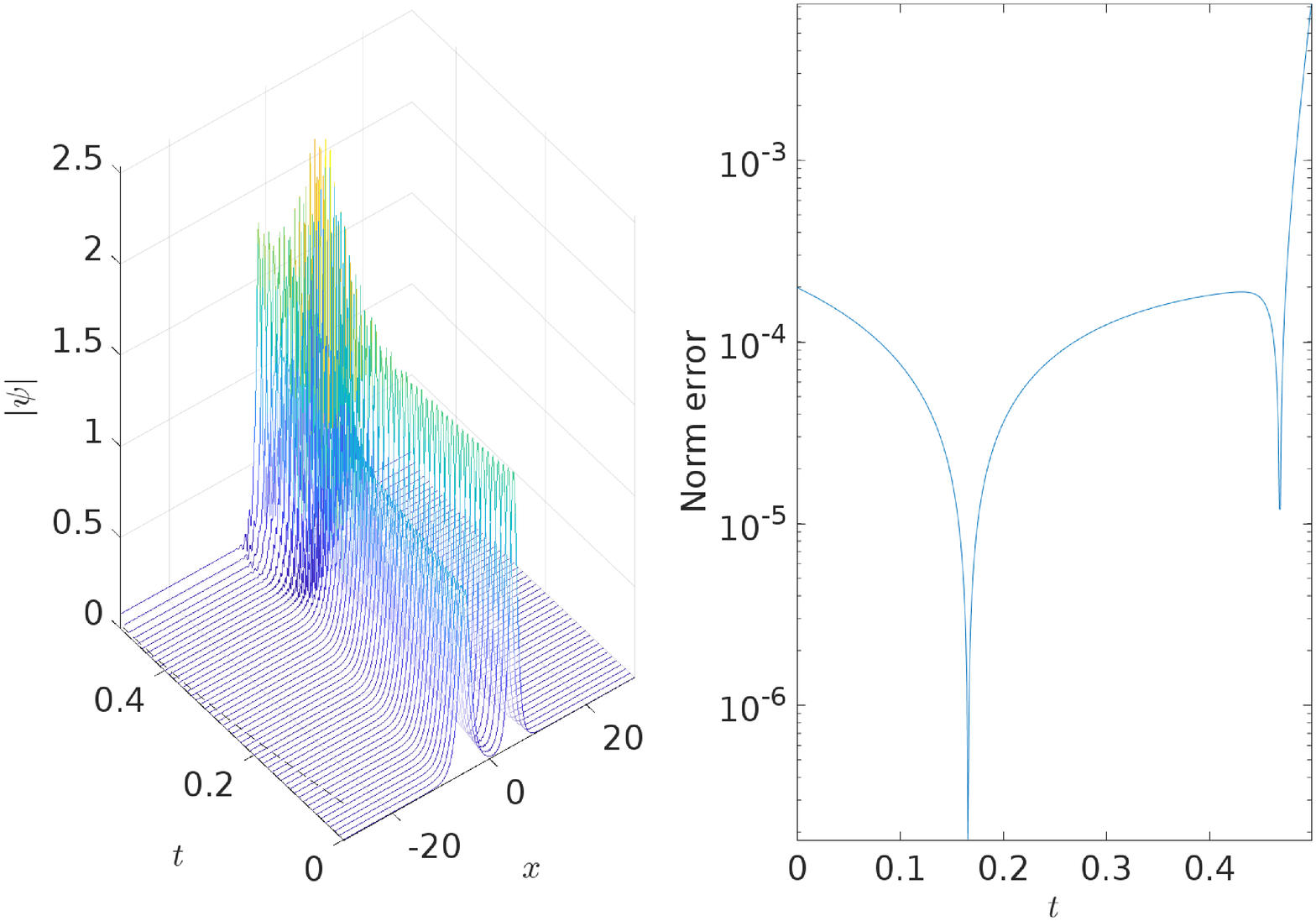}
 }
\caption{Left: Solution of (\ref{eq:nls5}) with $\psi(x,0) = \exp(8\ti x)\sech(x+5) +1.5\exp(-7\ti x)\sech(1.5(x-5))$ and 
periodic boundary conditions. Right: Error in norm.  (a) Results for the exponential method (\ref{eq:nlsdisc1})-(\ref{eq:nlsdisc2});  
(b) Results for the implicit midpoint time-discretization.}
 \label{fig:mixed_box_soliton_pair}
\end{figure}

Figure \ref{fig:mixed_box_soliton_pair} shows an example where the equation is initialized with propagating solitons to produce a collision. As a result, dispersive waves of much smaller amplitude emit out of the collision. 
In this experiment, the error in norm for the exponential time-discretization is compared to that of the underlying method 
(implicit midpoint in time, symplectic Euler in space).   Recent work proposed the implicit midpoint time-discretization for a similar 
equation \cite{HDY17}, showing that the discretization is almost structure-preserving, 
though it does preserve the norm in the absence of weak damping.  Only the exponential integrator shows the correct behavior when $\alpha$ and $\beta$ are nonzero.
%%and hence puts an stringent time step-size requirement to converge for longer simulations.

%%%%%%%%%%%%%%%%%%%%%%%%%%%%%%%%%%%%%%%%%%%%%%%%%%%%%%%%%%%%%%%%%%%%%%%%%%%%%%%%%%%%%%%%%%%%%%%%%%
\subsection{Damped-driven Camassa-Holm equation}
Consider the discretization (\ref{eq:expbox}) applied to equation (\ref{eq:dch}).
% A conformal multi-symplectic numerical method for PDE (\ref{eq:cmsdch}) is
% $${\bf K}(\prescript{\gamma/2}{}D_t^{\gamma/2} A_x z) + {\bf L} (D_x \prescript{\gamma/2}{}A_t^{\gamma/2} z) = \nabla S_\gamma (A_x \prescript{\gamma/2}{}A_t^{\gamma/2} z)$$
% This method can be reduced to
% \begin{align*}
% &\prescript{\gamma/2}{}D_t^{\gamma/2} A_{x^3} \prescript{\gamma/2}{}A_t^{\gamma/2} u -  \prescript{\gamma/2}{}D_t^{\gamma/2} D_{x^2} A_x \prescript{\gamma/2}{}A_t^{\gamma/2} u + 3 A_x \prescript{\gamma/2}{}A_t^{\gamma/2} (D_x A_x \prescript{\gamma/2}{}A_t^{\gamma/2} u ~~A_{x^2} \prescript{\gamma/2}{}A_t^{\gamma/2} u) \\
% &= 3 \prescript{\gamma/2}{}A_t^{\gamma/2} (D_x A_{x^2} \prescript{\gamma/2}{}A_t^{\gamma/2} u ~~D_{x^2} A_x \prescript{\gamma/2}{}A_t^{\gamma/2} u) - A_x \prescript{\gamma/2}{}A_t^{\gamma/2} (D_x A_x \prescript{\gamma/2}{}A_t^{\gamma/2} u ~~D_{x^2} \prescript{\gamma/2}{}A_t^{\gamma/2} u) \\
% &+ \prescript{\gamma/2}{}A_t^{\gamma/2} (A_{x^3} \prescript{\gamma/2}{}A_t^{\gamma/2}u ~~D_{x^3} \prescript{\gamma/2}{}A_t^{\gamma/2} u)
% \end{align*}
Take $A_{x^n} = A_x A_x A_x \ldots$ ($n$ times) and $\delta_{x}^{2} = \delta_{x}^{+} \delta_{x}^{-}$, 
and apply the method directly to (\ref{eq:dch}) to get
\begin{gather}
\begin{aligned} \label{eq:cimp-dch}
&D_{t_i}^{\gamma} A_{x^2} u_{n+1/2}^{i} - D_{t_i}^{\gamma} \delta_{x}^{2} u_{n+1/2}^{i} + 3A_x(\delta_{x}^{+} A_{t_i}^{\gamma}u_{n+1/2}^{i}~~ A_{x} A_{t_i}^{\gamma} u_{n+1/2}) \\
&= 3 (\delta_{x}^{+} A_{x} A_{t_i}^{\gamma} u_{n+1/2}^{i} ~~\delta_{x}^{2}  A_{t_i}^{\gamma} u_{n+1/2}^{i}) - A_x (\delta_{x}^{+} A_{t_i}^{\gamma} u_{n+1/2}^{i} ~~\delta_{x}^{2} A_{t_i}^{\gamma} u_{n}^{i}) 
+ (A_{x^2} A_{t_i}^{\gamma}u_{n+1/2}^{i} ~~\delta_{x}^{2} \delta_{x}^{+} A_{t_i}^{\gamma} u_{n}^{i}).
\end{aligned}
\end{gather}
It can be shown that this method preserves the Casimir $\int u \, dx \approx \sum_n u_n^i$.

\begin{figure}[htbp]
\centering
 \begin{tabular}{lll}
\includegraphics[width=0.5\textwidth]{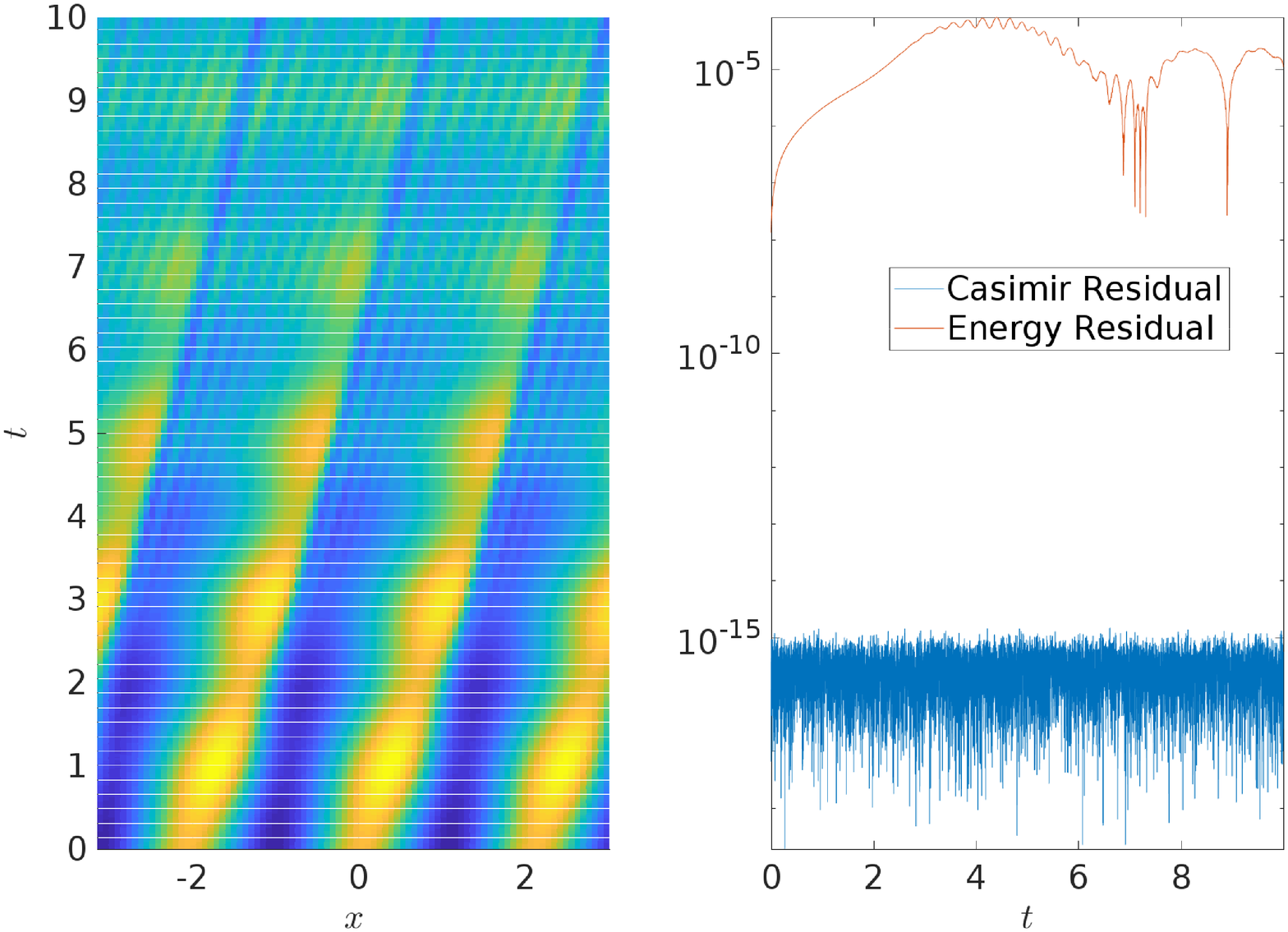} & 
\includegraphics[width=0.5\textwidth]{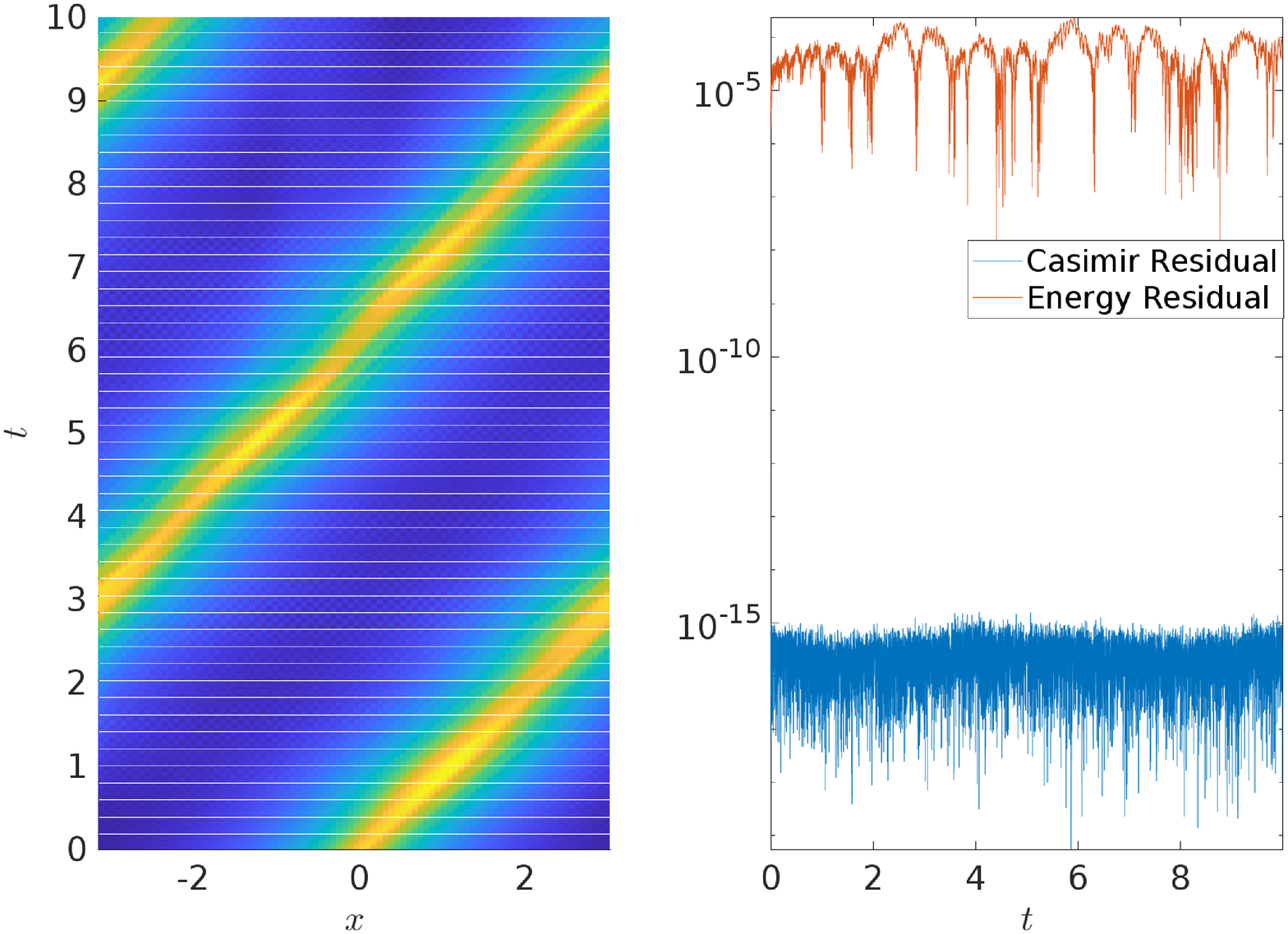} \\ 
\includegraphics[width=0.5\textwidth]{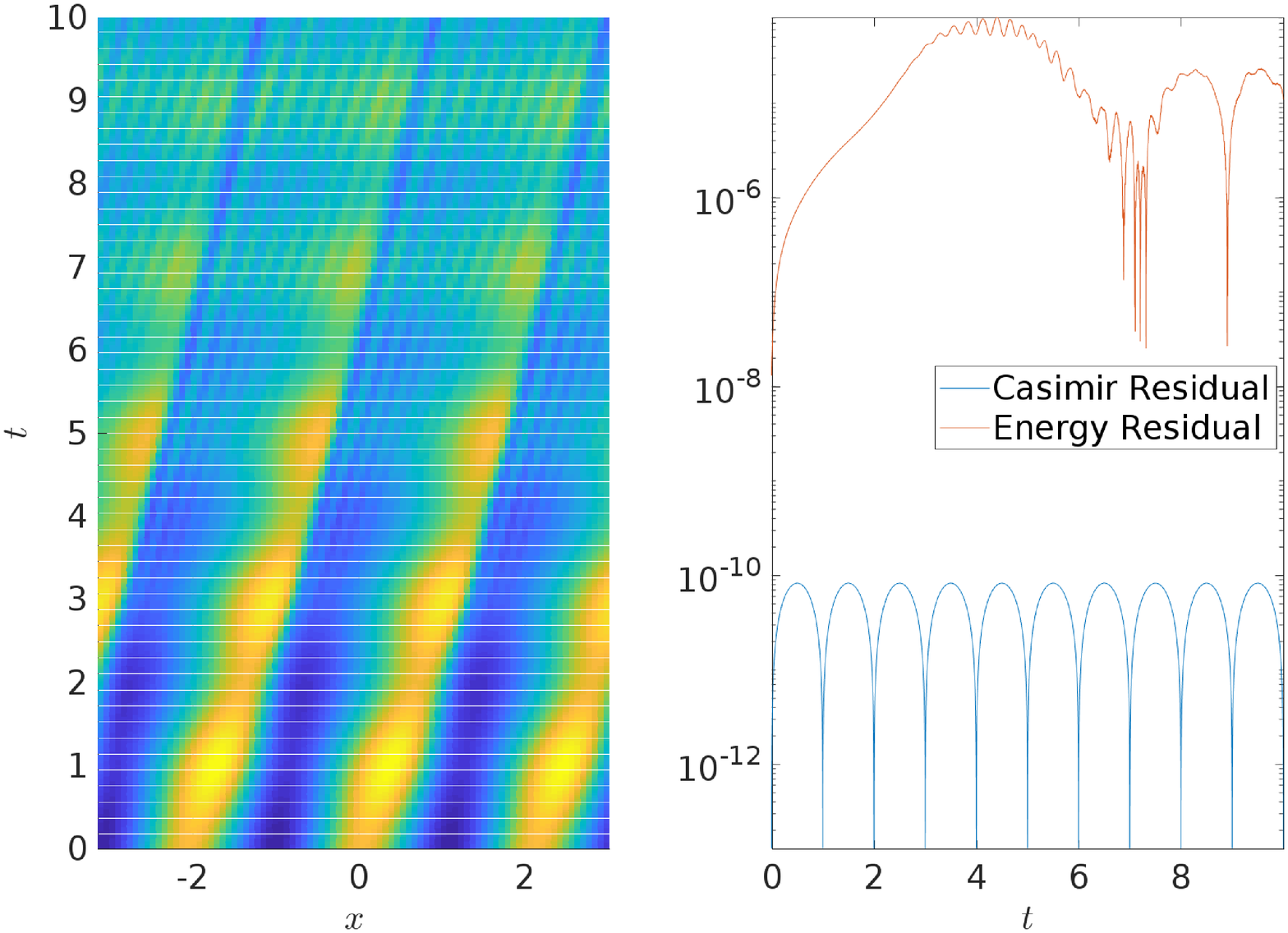} & 
\includegraphics[width=0.5\textwidth]{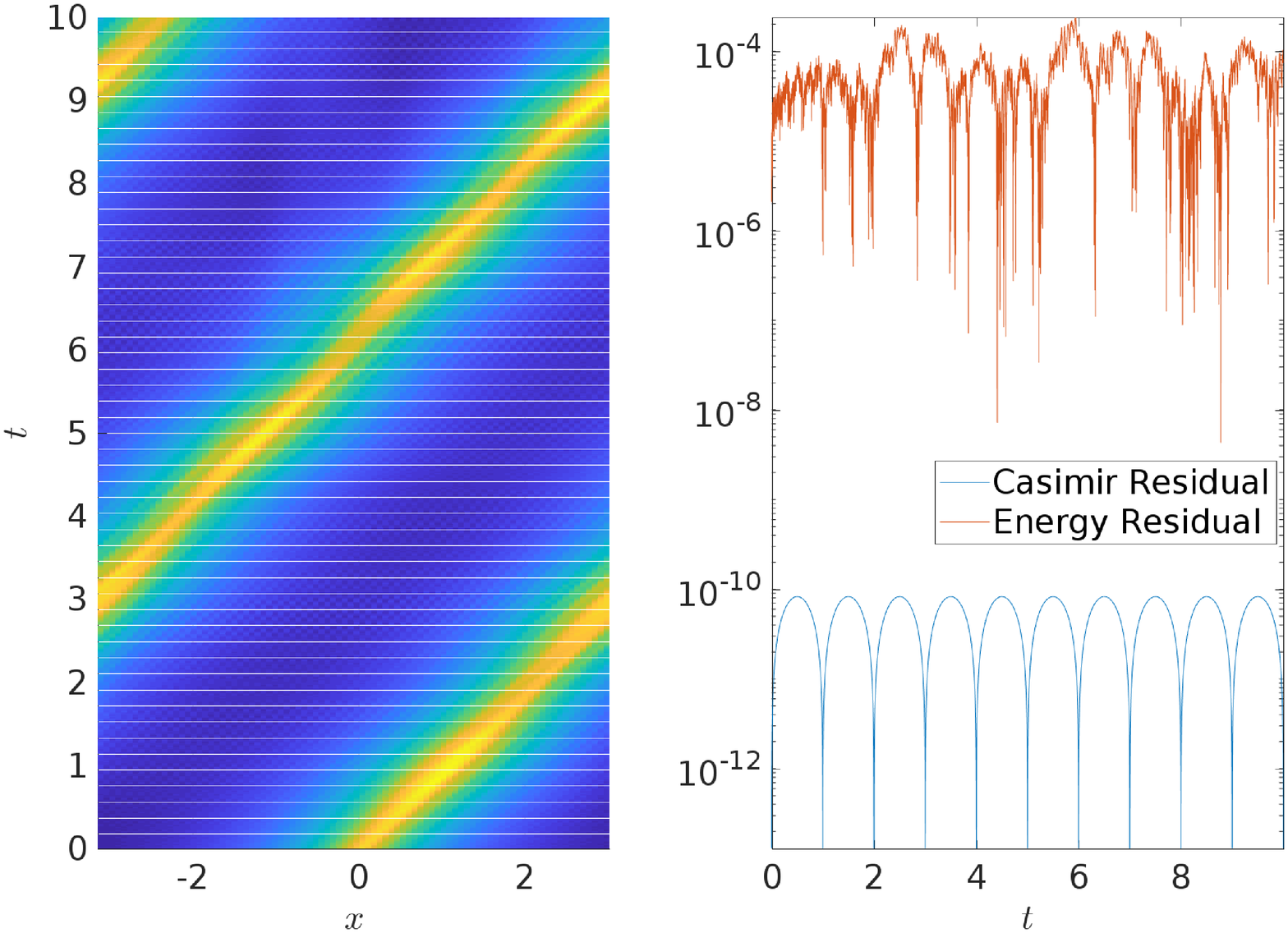}
 \end{tabular}
 \caption{Solutions of (\ref{eq:dch}) next to the Casimir and energy residuals.  
 The top row shows results from the conformal multi-symplectic method (\ref{eq:cimp-dch}), 
 while the bottom row shows results from the underlying method, the Preissmann box scheme.  
 Left columns show results with $u_0  = 0.2 + 0.1\cos(3x)$, 
 and the right columns show results for $u_0 = e^{-|x|}$.}
 \label{fig:Casimir-and-Energy-Residual}
 \end{figure}

Figure~\ref{fig:Casimir-and-Energy-Residual}, shows the residual in the Casimir and the Energy over the time interval $[0,10]$
and spatial domain $[-\pi,\pi]$ when method (\ref{eq:cimp-dch}) is used to propagate smooth and non-smooth initial data
\begin{align} \label{eq:ic}
u_0  = 0.2 + 0.1\cos(3x), ~\text{and}~ u_0 = e^{-|x|}
\end{align}
with
\begin{align*}
\dx = 0.07, ~\dt =0.001, ~x \in [-\pi,\pi], ~t \in [0,10], ~\gamma(t) = -0.2\sin(\pi t), ~\text{and}~f(x,t)=0.
\end{align*}
The figures confirm Casimir preservation by the structure-preserving method (\ref{eq:cimp-dch}). 
As expected, the non-smooth kink initial profile $e^{-|x|}$ evolves and produces persisting unstable modes early in the simulation 
because of the non-preservation of the decay in the energy or the $\IH^1$ norm of the solution. 
The non-smooth solution maintains its shape and direction throughout the numerical simulations nonetheless. 
Although neither method preserves the energy of the system, the exponential integrator 
(\ref{eq:cimp-dch}) has energy residual one order of magnitude lower than the underlying method. 
The result of the periodic input to the system is evident in all the plots. 
The solutions' amplitude and the energy residual due to the underlying method varies periodically in sync with the amplitude of $\gamma$.

%%%%%%%%%%%%%%%%%%%%%%%%%%%%%%%%%%%%%%%%%%%%%%%%%%%%%%%%%%%%%%%%%%%%%%%%
\section{Conclusion} \label{sec:conc}
%%%%%%%%%%%%%%%%%%%%%%%%%%%%%%%%%%%%%%%%%%%%%%%%%%%%%%%%%%%%%%%%%
We have presented general systems of PDEs with time-dependent coefficients, along with damping/driving forces,
that have solutions which satisfy certain conservation laws.  We have shown that several PDEs arising from physical applications have 
the general form we consider, including damped/driven KdV, nonlinear wave, NLS, and Camassa-Holm equations.
Since these equations have conservation laws that are desirable to preserve numerically, we present three, second-order in time, 
discretizations that preserve certain conservation laws.  Even in cases where the damping coefficients are constant, 
the preservation of these local conservation laws signifies a stronger result that the global (boundary-condition-dependent) 
conservation that was previously achieved with conformal symplectic methods.  Ultimately, the number of PDEs that might be 
considered under this framework, along with the number of structure-preserving discretizations we might propose, is vast.  
To demonstrate the effectiveness of our approach, we numerically solve a damped-driven NLS equation and a damped-driven Camassa-Holm equation 
using discretizations that preserve particular properties of the equations.  In addition to illustrating preservation of the
desired properties, the results exhibit qualitatively correct behavior in other respects, such as energy, as well as 
favorable comparison to other methods that have been proposed for similar equations.
 
%% References


\bibliographystyle{elsarticle-num}

\begin{thebibliography}{00}

%\bibitem{ARZ13} M.~Asadzadeh, D.~Rostamy, F.~Zabihi, Discontinuous Galerkin and multiscale variational schemes for a coupled damped nonlinear system of Schr\"{o}dinger equations, {\sl Numer.~Meth.~Part.~D.~E.} 29 (2013) 1912-1945.

\bibitem{BZ04} I.V.~Barashenkov and E.V.~Zemlyanaya, Traveling solitons in the damped-driven nonlinear Schr\"{o}dinger equation, {\sl SIAM J.~Appl.~Math.} 64 (2004) 800-818.

%\bibitem{B16} A.~Bhatt, Multi-conformal-symplectic numerical methods for damped driven nonlinear Schr\"{o}dinger equations, {\sl submitted}, 2016.

\bibitem{BFM15} A.~Bhatt, D.~Floyd, and B.E.~Moore, Second order conformal symplectic schemes for damped Hamiltonian systems, {\sl J.~Sci.~Comput.} 66 (2016) 1234-1259.

\bibitem{BM17} A.~Bhatt and B.E.~Moore, Structure-preserving exponential Runge-Kutta methods, {\sl SIAM J.~Sci.~Comput.} to appear, 2017.

\bibitem{B09} A.~Biswas, Solitary wave solution for the generalized KdV equation with time-dependent damping and dispersion, {\sl Commun.~Nonlinear Sci.}, 14 (2009) 3503-3506.

%\bibitem{BS94} P.B.~Bochev and C.~Scovel, On quadratic invariants and symplectic structure, {\sl BIT} 34.3 (1994) 337-345.

\bibitem{CZW17} W.~Cai, H.~Zhang, and Y.~Wang, Modelling damped acoustic waves by a dissipation-preserving conformal symplectic method, {\sl Proc.R.Soc.A} 473:20160798 (2017).

\bibitem{Cetal12} E.~Celledoni, V.~Grimm, R.I.~McLachlan, D.I.~McLaren, D.~O'Neale, B.~Owren, G.R.W.~Quispel, 
Preserving energy resp. dissipation in numerical PDEs using the “Average Vector Field” method, {\sl J.~Comp.~Phys.} 231.20 (2012) 6770-6789.

\bibitem{CO08}
D.~Cohen, B.~Owren, and X.~Raynaud. Multi-symplectic integration of the Camassa–Holm equation. {\sl Journal of Computational Physics} 227(2008) 5492-5512.

%\bibitem{CMMOQW} E.~Celledoni, R.I.~McLachlan, D.I.~McLaren, B.~Owren, G.R.W.~Quispel, and W.M.~Wright, ``Energy-preserving Runge-Kutta Methods,'' 
%{\sl ESAIM: Mathematical Modeling and Numerical Analysis}, 43:645-649, 2009.
%
%\bibitem{CCO08} Celledoni, Elena, David Cohen, and Brynjulf Owren. 
%``Symmetric exponential integrators with an application to the cubic Schr{\"o}dinger equation." Foundations of Computational Mathematics 8.3 (2008): 303-317.

\bibitem{DA15} M.~D'Abbicco, The threshold of effective damping for semilinear wave equations, {\sl Math.~Methods Appl.~Sci.} 38 (2015) 1032–1045.

\bibitem{DALR13} M.~D'Abbicco, S.~Lucente, M.~Reissig, Semi-linear wave equations with effective damping, {\sl Chin.~Ann.~Math.~Ser.~B} 34 (2013) 345–380.

\bibitem{FQZ05} Y.~Feng, W.-X.~Qin, and Z.~Zheng, Existence of localized solutions in the parametrically driven and damped DNLS equation in high dimensional lattices, 
{\sl Phys.~Lett.~A}, 346 (2005) 99-110.

%\bibitem{FS97} W.I.~Fushchych, Z.I.~Symenoh, Symmetry of equations with convection terms, {\sl J.~Nonlinear Math.~Phy.} 4 (1997) 470-479.

\bibitem{Fu16} H.~Fu, W.-E.~Zhou, X.~Qian, S.-H.~Song, and L.-Y.~Zhang, Conformal structure-preserving method for damped nonlinear Schr\"{o}dinger equation, {\sl Chin.~Phys.~B} 25.11 (2016) 110201.

\bibitem{HR06}
H.~Holden, X.~Raynaud. A convergent numerical scheme for the Camassa–Holm equation based on multipeakons. {\sl Discrete and Continuous Dynamical Systems} 14 (2006) 505–523.

\bibitem{HDY17} W.~Hu, Z.~Deng, and T.~Yin, Almost structure-preserving analysis for weakly linear damping nonlinear Schrödinger equation with periodic perturbation,
{\sl Commun.~Nonlinear Sci.~Numer.~Simulat.} 42 (2017) 298-312. 

\bibitem{ChiPhys02} P.~Jie-Hua, T.~Jia-Shi, Y.~De-Jie, Y.~Jia-Ren, and H.~Wen-Hua, Soluitons, bifurcations and chaos of the nonlinear Schro\"{o}dinger equation with weak damping,
{\sl Chin.~Phys.} 11 (2002) 213-217.

\bibitem{JK10} A.G.~Johnpillai and C.M.~Khalique, Group analysis of KdV equation with time dependent coefficients, {\sl App.~Math.~Comput.} 216 (2010) 3761–3771

%%\bibitem{Kl08} Klein, Christian. ``Fourth order time-stepping for low dispersion Korteweg-de Vries and nonlinear Schr{\"o}dinger equation." Electronic Transactions on Numerical Analysis 29 (2008): 116-135.
\bibitem{KL98}
Y.~ S.~Kivshara, and B.~ Luther-Daviesb. Dark optical solitons: physics and applications. {\sl Physics reports} 298.2-3 (1998) 81-197.

\bibitem{KS00}
S.~Kouranbaeva, S.~Shkoller. A variational approach to second-order multisymplectic field theory. {\sl Journal of Geometry and Physics} 35 (2000) 333-366.

\bibitem{KPH} V.I.~Kruglov, A.C.~Peacock, J.D.~Harvey, Exact solutions of the generalized nonlinear Schrödinger equation with distributed coefficients, {\sl Phys.~Rev.~E Stat.~Nonlin.~Soft Matter Phys.} 71(5 Pt 2):056619, 2005.

\bibitem{KS01} J.G.~Kingstona, C.~Sophocleousb, Symmetries and form-preserving transformations of one-dimensional wave equations with dissipation, {\sl Int.~J.~Nonlin.~Mech.} 36 (2001) 987–997.

\bibitem{Law67} D.J.~Lawson, Generalized Runge-Kutta processes for stable systems with large Lipschitz constants. {\sl SIAM J.~Numer.~Anal.} 4.3 (1967) 372-380.

\bibitem{MP01} R.I.~McLachlan and M.~Perlmutter. Conformal Hamiltonian systems. {\sl J.~Geom.~Phys.} 39.4 (2001) 276-300.

\bibitem{McLQ01} R.I.~McLachlan, G.R.W.~Quispel,  What kinds of dynamics are there? Lie pseudogroups, dynamical systems and geometric integration, {\sl Nonlinearity} 14:1689--1705 (2001).

\bibitem{McLQ02} R.I.~McLachlan, G.R.W.~Quispel, Splitting methods, {\sl Acta Numer.} 11:341--434 (2002).
%
%\bibitem{McLRS14} R.I.~McLachlan, B.N.~Ryland, and Y.~Sun, ``High-order multisymplectic Runge-Kutta methods,'' {\sl SIAM J.~Sci.~Comput.} 36(5):A2199-A2226, 2014

\bibitem{MS11} K.~Modin and G.~S{\"o}derlind, Geometric integration of Hamiltonian systems perturbed by Rayleigh damping, {\sl BIT} 51.4 (2011): 977-1007.

\bibitem{thesis} B.E. Moore. {\em A Modified Equations Approach for Multi-Symplectic Integration Methods}. PhD Thesis, University of Surrey, 2003.

\bibitem{Moore09} B.E.~Moore. Conformal multi-symplectic integration methods for forced-damped semi-linear wave equations, {\sl Math.~Comput.~Simulat.} 80 (2009) 20-28.

\bibitem{M16} B.E.~Moore, Multi-conformal-symplectic PDEs and discretizations, {\sl J.~Comput.~Appl.~Math.}, 323:1-15, 2017..

\bibitem{MNS13} B.E.~Moore, L.~Noreña, and C.M.~Schober, Conformal conservation laws and geometric integration for damped Hamiltonian PDEs, {\sl J.~Comp.~Phys.} 232 (2013) 214-233.

\bibitem{N11} K.~Nishihara, Asymptotic behavior of solutions to the semilinear wave equation with time-dependent damping, {\sl Tokyo J.~Math.}, 34 (2011) 327-343.

\bibitem{NZ09} K.~Nishihara and J.~Zhai, Asymptotic behaviors of solutions for time dependent damped wave equations, {\sl J.~Math.~Anal.~Appl.} 360 (2009) 412–421.
%
%\bibitem{R00} S.~Reich, ``Multi-symplectic Runge-Kutta collocation methods for Hamiltonian wave equations,'' {\sl J.~Comp.~Phys.}, 157:473-499, 2000.

%\bibitem{S05} H.~Segur, D.~Henderson, J.~Carter, J.~Hammack, C.-M.~Li, D.~Pheiff, K.~Socha, Stabilizing the Benjamin-Feir instability, {\sl J.~Fluid Mech.} 539 (2005) 229-271.

\bibitem{SH} V.N.~Serkin, A.~Hasegawa, Novel soliton solutions of the nonlinear Schrodinger equation model {\sl Phys.~Rev.~Lett.} 85(21):4502-5, 2000.

\bibitem{SMB} V.N.~Serkin, M.~Matsumoto, T.L.~Belyaeva, Bright and dark solitary nonlinear Bloch waves in dispersion managed fiber systems and soliton lasers,
{\sl Optics Communications}, 19 (2001) 159-171.

\bibitem{Shen} C.~Shen, Time-periodic solution of the weakly dissipative Camassa-Holm equation, 
{\sl Inter.~J.~Differential Equ.} 2011 (2011), Article ID 463416, 16 pages

\bibitem{SS06}
M.~Stanislavova, and A.~Stefanov. Attractors for the viscous Camassa-Holm equation. {\sl arXiv preprint math/0612321} 2006.

\bibitem{SQWS10} H.~Su, M.~Qin, Y.~Wang, R.~Scherer, Multi-symplectic Birkhoffian structure for PDEs with dissipation terms, {\sl Phys.~Lett.~A}, 374 (2010) 2410-2416.

\bibitem{SS05} Y.~Sun and Z.~Shang, Structure-preserving algorithms for Birkhoffian systems, {\sl Phys.~Lett.~A}, 336 (2005) 358-369.

\bibitem{SB02} V.S.~Shchesnovich and I.V.~Barashenkov, Soliton-radiation coupling in the parametrically driven, damped nonlinear Schr\"{o}dinger equation, {\sl Physica D} 164 (2002) 83-109.

%\bibitem{Burgh02} van der Burgh, A. H. P. An equation with a time-periodic damping coefficient: stability diagram and an application. Delft University of Technology, Faculty of Electrical Engineering, Mathematics and Computer Science, Delft Institute of Applied Mathematics, 2002.

\bibitem{W17} Y.~Wakasugi, Scaling variables and asymptotic profiles for the semilinear damped wave equation with variable coefficients, {\sl J.~Math.~Anal.~Appl.} 447 (2017) 452–487.

\bibitem{W04} J.~Wirth, Solution representations for a wave equation with weak dissipation, {\sl Math.~Meth.~Appl.~Sci.} 27 (2004) 101-124.

\bibitem{WY06}
S.~Wu, and Z.~Yin. Blow-up, blow-up rate and decay of the solution of the weakly dissipative Camassa–Holm equation. {\sl J. Math. Phys.} 47(2006) 1–12.

\bibitem{WY07}
S.~Wu, and Z.~Yin. Blowup and decay of solution to the weakly dissipative Camassa–Holm equation. {\sl Acta Math. Appl. Sin.} 30 (2007) 996–1003.

\bibitem{WY09}
S.~Wu, and Z.~Yin. Global existence and blow-up phenomena for the weakly dissipative Camassa–Holm equation. {\sl J. Differential Equations} 246 (2009) 4309–4321.

\bibitem{XFS06} T.~Xiao-Yan, H.~Fei, and L.~Sen-Yue, Variable coefficient KdV equation and the analytical diagnoses of a dipole blocking life cycle, {\sl Chinese Phys.~Lett.} 23 (2006) 887-890.
 \end{thebibliography}



\end{document}